\documentclass[11pt, twoside]{amsart}

\usepackage{t1enc}
\usepackage[latin1]{inputenc}
\usepackage{latexsym}
\usepackage{amssymb}
\usepackage{graphicx}
\usepackage{xcolor}
\usepackage{amsmath}
\usepackage{amsthm}
\usepackage{amsfonts}
\usepackage{mathpazo}
\usepackage{mathrsfs}
\usepackage[paper = letterpaper, left = 3.5cm, right = 3.5cm, headsep = 6mm, footskip = 10mm,
top = 35mm, bottom = 35mm, footnotesep=5mm, headheight = 2cm]{geometry}
\usepackage{fancyhdr}
\usepackage[all]{xy}
\usepackage[british]{babel}
\usepackage{soul}
\usepackage{hyperref}

\newcommand*{\widebar}{\overline}
\newcommand*{\definiere}{\mathrel{\mathop:}=}

\newcommand*{\tensor}{\otimes}

\newcommand{\cyclic}{\mathop{\kern0.9ex{{+}\kern-2.10ex\raise-0.20
      ex\hbox{\Large\hbox{$\circlearrowright$}}}}\limits}
\newcommand{\acts}{\mbox{ \raisebox{0.26ex}{\tiny{$\bullet$}} }}

\def\N{\ifmmode{\mathbb N}\else{$\mathbb N$}\fi}
\def\Z{\ifmmode{\mathbb Z}\else{$\mathbb Z$}\fi}
\def\Q{\ifmmode{\mathbb Q}\else{$\mathbb Q$}\fi}
\def\R{\ifmmode{\mathbb R}\else{$\mathbb R$}\fi}
\def\C{\ifmmode{\mathbb C}\else{$\mathbb C$}\fi}
\def\K{\ifmmode{\mathbb K}\else{$\mathbb K$}\fi}
\def\P{\ifmmode{\mathbb P}\else{$\mathbb P$}\fi}
\def\g{\ifmmode{\mathfrak g}\else {$\mathfrak g$}\fi}
\def\h{\ifmmode{\mathfrak h}\else {$\mathfrak h$}\fi}
\def\a{\ifmmode{\mathfrak a}\else {$\mathfrak a$}\fi}
\def\k{\ifmmode{\mathfrak k}\else {$\mathfrak k$}\fi}
\def\p{\ifmmode{\mathfrak p}\else {$\mathfrak p$}\fi}
\def\b{\ifmmode{\mathfrak b}\else {$\mathfrak b$}\fi}
\def\n{\ifmmode{\mathfrak n}\else {$\mathfrak n$}\fi}
\def\m{\ifmmode{\mathfrak m}\else {$\mathfrak m$}\fi}
\def\t{\ifmmode{\mathfrak t}\else {$\mathfrak t$}\fi}
\def\O{\ifmmode{\mathscr{O}}\else {$\mathscr{O}$}\fi}
\def\W{\ifmmode{\mathcal{V}}\else {$\mathscr{W}$}\fi}
\def\id{{\rm id}}

\def\hq{/\hspace{-0.14cm}/}
\def\kleinematrix#1,#2,#3,#4,{\begin{pmatrix}#1 & #2 \\ #3 & #4
  \end{pmatrix}}

\DeclareMathOperator{\Int}{Int}
\DeclareMathOperator{\Ad}{Ad}

\DeclareMathOperator{\Lie}{Lie}

\DeclareMathOperator{\dom}{dom}

\DeclareMathOperator{\Type}{Type}

\newtheoremstyle{daniel}{3.0mm}{0mm}{\itshape}{}{\bfseries}{.}{1.5mm}{}
\theoremstyle{daniel}
\newtheorem{thm}{Theorem}[section]
\newtheorem{prop}[thm]{Proposition}
\newtheorem{Defi}[thm]{Definition}
\newtheorem{lemma}[thm]{Lemma}
\newtheorem{cor}[thm]{Corollary}
\newtheorem{Exs}[thm]{Examples}
\newtheorem{Ex}[thm]{Example}
\newtheorem{Rems}[thm]{Remarks}
\newtheorem{Rem}[thm]{Remark}

\newtheorem*{thm*}{Theorem}
\newtheorem*{cor*}{Corollary}

\newtheorem*{thm1.1}{Theorem 1.1}
\newtheorem*{thm11.5}{Theorem 11.5}
\newtheorem*{cor11.4}{Corollary 11.4}
\newtheorem*{cor11.6}{Corollary 11.6}
\newtheorem*{prop*}{Proposition}
\newtheorem*{claim}{Claim}

\newenvironment{rem}   {\begin{Rem}\em}{\end{Rem}}
\newenvironment{rems}   {\begin{Rems}\em}{\end{Rems}}

\newenvironment{ex}  {\begin{Ex}\em}{\end{Ex}}

\newenvironment{enumerate-p}{\begin{enumerate}}{
 \setcounter{equation}{\value{enumi}}
\end{enumerate}}

\def\cprime{$'$} \def\polhk#1{\setbox0=\hbox{#1}{\ooalign{\hidewidth
  \lower1.5ex\hbox{`}\hidewidth\crcr\unhbox0}}}
  \def\polhk#1{\setbox0=\hbox{#1}{\ooalign{\hidewidth
  \lower1.5ex\hbox{`}\hidewidth\crcr\unhbox0}}}
\providecommand{\bysame}{\leavevmode\hbox to3em{\hrulefill}\thinspace}
\providecommand{\MR}{\relax\ifhmode\unskip\space\fi MR }
\providecommand{\MRhref}[2]{%
  \href{http://www.ams.org/mathscinet-getitem?mr=#1}{#2}
}
\providecommand{\href}[2]{#2}
\begin{document}
\title[Compact Kaehler quotients and Geometric Invariant Theory]{Compact Kaehler quotients of algebraic varieties and Geometric Invariant Theory}
\author{Daniel Greb}
\date{}
\address{Mathematisches Institut\\
Abteilung f\"ur Reine Mathematik\\
Albert-Ludwigs-Universit\"at\\
Eckerstr. 1\\
79104 Freiburg im Breisgau\\
Germany}
\email{daniel.greb@math.uni-freiburg.de}{}
\urladdr{\href{http://home.mathematik.uni-freiburg.de/dgreb}{http://home.mathematik.uni-freiburg.de/dgreb}}

\begin{abstract}
Given an action of a complex reductive Lie group $G$ on a normal variety $X$, we show that every analytically Zariski-open subset of $X$ admitting an analytic Hilbert quotient with projective quotient space is given as the set of semistable points with respect to some $G$-linearised Weil divisor on $X$. Applying this result to Hamiltonian actions on algebraic varieties, we prove that semistability with respect to a momentum map is equivalent to GIT-semistability in the sense of Mumford and Hausen. It follows that the number of compact momentum map quotients of a given algebraic Hamiltonian $G$-variety is finite. As further corollary we derive a projectivity criterion for varieties with compact Kaehler quotient.
\end{abstract}
\maketitle

\section{Introduction and statement of main results}
Based on the fundamental work of Guillemin-Sternberg~\cite{GuilleminGeometricQuantisation}, Kirwan~\cite{KirwanCohomology}, Mumford, and others, momentum geometry has become one of the most important tools for studying actions of complex-reductive Lie groups $G=K^\C$ on complex spaces. Given a K\"ahlerian holomorphic $G$-space and a momentum map $\mu: X \to \Lie(K)^*$ with respect to a $K$-invariant K\"ahler form $\omega$, the set of \emph{$\mu$-semistable points} $X(\mu) \definiere \{x \in X \mid \overline{G\acts x} \cap \mu^{-1}(0) \neq \emptyset \}$ admits an analytic Hilbert quotient, i.e., a $G$-invariant holomorphic Stein map $\pi: X(\mu)\to
X(\mu)\hq G$ onto a K\"ahlerian complex space $X(\mu)\hq G$ with structure sheaf $\mathscr{O}^{hol}_{X(\mu)/\negthickspace / G}=(\pi_*\mathscr{O}_{X(\mu)}^{hol})^G$, see \cite{ReductionOfHamiltonianSpaces}, \cite{Extensionofsymplectic}, and \cite{SjamaarSlices}. If $X$ is projective algebraic, and if the K\"ahler form $\omega$ as well as the momentum map $\mu$ are induced by an embedding of $X$ into some projective space, both the set $X(\mu)$ of semistable points and the quotient $X(\mu)\hq G$ can also be constructed via Geometric Invariant Theory (GIT). In particular, the complex space $X(\mu)\hq G$ is projective algebraic, and the map $\pi: X(\mu) \to X(\mu)\hq G$ is a good quotient in the sense of GIT, cf.\ \cite{MumfordGIT}.

However, already on projective manifolds there exist many K\"ahler forms which do not arise as curvature forms of ample line bundles. From the point of view of complex and symplectic geometry it is therefore natural to study also the semistability conditions induced by these forms and their relation to the algebraic geometry of the underlying $G$-variety. Furthermore, especially when studying questions related to the variation of GIT quotients and applications, e.g.\ to moduli problems, one encounters interesting phenomena related to non-integral K\"ahler forms, even if one is a priori interested in ample line bundles. See \cite{SchmittGieseker} for an example related to moduli spaces of semistable sheaves on higher-dimensional projective manifolds. The above discussion motivates the following definition: an \emph{algebraic Hamiltonian $G$-variety} is a complex algebraic $G$-variety $X$ together with a (not necessarily integral) $K$-invariant K\"ahler form $\omega$ and a $K$-equivariant momentum map $\mu: X \to \Lie(K)^*$ with respect to $\omega$.

Under certain mild assumptions on the singularities, it was shown in \cite{PaHq} that compact momentum map quotients of algebraic Hamiltonian $G$-varieties are projective algebraic and that the corresponding sets of semistable points are algebraically Zariski-open. This already is quite remarkable in view of examples of compact non-projective geometric quotients of smooth projective algebraic $\C^*$-varieties constructed by Bia{\l}ynicki-Birula and {\'S}wi{\polhk{e}}cicka \cite{BBExoticOrbitSpaces}.

In this paper we show that momentum map quotients of algebraic Hamiltonian $G$-varieties have even stronger algebraicity properties and we give an essentially complete picture of the relation between momentum geometry and Geometric Invariant Theory for Hamiltonian actions on algebraic varieties. We conclude this introduction with a summary of the main results and an outline of the paper.

\begin{thm}[Algebraicity of momentum map quotients]\label{thm:mainthmHamiltonian}
Let $G=K^\C$ be a complex reductive Lie group and let $X$ be a $G$-irreducible algebraic Hamiltonian $G$-variety with at worst $1$-rational singularities. Assume that the zero fibre $\mu^{-1}(0)$ of the momentum map $\mu: X \to \Lie(K)^*$ is nonempty and compact. Then
\begin{enumerate}
\item the analytic Hilbert quotient $X(\mu)\hq G$ is a projective algebraic variety,
\item the set $X(\mu)$ of $\mu$-semistable points is algebraically Zariski-open in $X$,
\item the map $\pi: X(\mu) \to X(\mu)\hq G$ is a good quotient,
\item there exists a $G$-linearised Weil divisor $D$ (in the sense of \cite{HausenGITwithWeildivisors}) such that $X(\mu)$ coincides with the set $X(D, G)$ of semistable points with respect to $D$.
\end{enumerate}
\end{thm}
Parts (1) and (2) are already contained in \cite{PaHq}, and we have included them here in order to convey a complete picture of the situation; the main contribution of this paper is part (3). Theorem~\ref{thm:mainthmHamiltonian} generalises results proven by Heinzner and Migliorini \cite{MomentumProjectivity} for smooth projective Hamiltonian $G$-varieties to a singular, non-projective and non-compact setup. Note that Theorem~\ref{thm:mainthmHamiltonian} applies in particular to Hamiltonian $G$-varieties with proper momentum map.

A variety or complex space $X$ is said to have only \emph{1-rational singularities}, if for any re\-so\-lu\-tion of singularities $f: \widetilde X \to X$ the sheaf $R^1f_*\mathscr{O}_{\widetilde X}$ vanishes. The class of complex spaces with $1$-rational singularities is the natural class of singular spaces to which projectivity results for K\"ahler Moishezon manifolds generalise, cf.\ \cite{ProjectivityofMoishezon}. Furthermore, it is stable under taking good quotients (in the algebraic category \cite{GrebSingularities}) and analytic Hilbert quotients (in the analytic category \cite{GrebAnalyticSingularities}).

In Section~\ref{sect:KaehlerquotientsGIT} we construct a K\"ahlerian non-projective proper algebraic surface, see Example~\ref{ex:nonprojective}, which shows the necessity of the assumption on the singularities in Theorem~\ref{thm:mainthmHamiltonian}. To the author's knowledge this surface is the first example of a K\"ahlerian non-projective proper algebraic variety in the literature.

Theorem~\ref{thm:mainthmHamiltonian} above is deduced from the following main result of this paper:
\begin{thm}\label{thm:mainthm}
Let $G$ be a complex reductive Lie group, let $X$ be a $G$-irreducible normal algebraic $G$-variety, and let $U\subset X$ be a $G$-invariant analytically Zariski-open subset of $X$ such that the analytic Hilbert quotient $\pi: U \to U\hq G$ exists. If $U\hq G$ is a projective algebraic variety, then
\begin{enumerate}
\item the set $U$ is algebraically Zariski-open in $X$,
\item the map $\pi: U \to U\hq G$ is is a good quotient.
\end{enumerate}
\end{thm}
In fact, large parts of the statement of Theorem~\ref{thm:mainthm} still hold under the weaker assumption that that the quotient $U\hq G$ is a complete algebraic variety as will become clear in the following outline of its proof.

As a first step we prove in Section~\ref{sect:Zariskiopen}:
\begin{thm*}[Openness Theorem]
Let $G$ be a connected complex reductive Lie group and let $X$ be an irreducible normal algebraic $G$-variety. Let $U$ be a $G$-invariant analytically Zariski-open subset of $X$ such that the analytic Hilbert quotient $\pi: U \to U \hq G$ exists. If $U\hq G$ is a complete algebraic variety, then $U$ is Zariski-open in $X$.
\end{thm*}
The main tool used in this section is the Rosenlicht quotient of $X$ by $G$, cf.\ Section~\ref{subsect:Rosenlicht}. A similar result has been proven by the author for Hamiltonian $G$-varieties and projective quotient spaces in \cite[Thm.\ 2]{PaHq}.

As a second step we prove algebraicity of the quotient map in Section~\ref{sect:algebraicmap}:
\begin{thm*}[Algebraicity Theorem]
Let $G$ be a connected complex reductive Lie group, let $X$ be an irreducible normal algebraic $G$-variety, and let $U$ be a $G$-invariant Zariski-open subset of $X$ such that the analytic Hilbert quotient $\pi: U \to Q$ exists. If $Q$ is a complete algebraic variety, then the quotient map $\pi$ is algebraic.
\end{thm*}
This result is new even for Hamiltonian actions and quotients of semistable points with respect to some momentum map. For the proof we construct an algebraic family of compact cycles from the group action and investigate the interplay between the associated map (to an auxiliary Rosenlicht quotient) and the quotient map $\pi$.

Section~\ref{sect:coherence} is devoted to the study of coherence properties of sheaves of invariants:
\begin{thm*}[Coherence Theorem]
Let $G$ be a connected complex reductive Lie group. Let $X$ be an irreducible algebraic $G$-variety with analytic Hilbert quotient $\pi: X \to Q \definiere X\hq G$ where $Q$ is a complete algebraic variety. If $\mathscr{F}$ is any coherent algebraic $G$-sheaf on $X$, then $(\pi_* \mathscr{F})^G$ is a coherent algebraic sheaf on $Q$.
\end{thm*}
The most important tools used here are Proposition~\ref{prop:meromorphicdescent} and Corollary~\ref{cor:meromorphicdescentgenericcase} concerning the descent of invariant meromorphic functions to analytic Hilbert quotients of not necessarily reduced holomorphic $G$-spaces.

Finally, the proof of Theorem~\ref{thm:mainthm} is given in Section~\ref{sect:maintheoremproof}. Here, we provide a short sketch. Assuming projectivity of the quotient $U\hq G$, Theorem~3 of \cite{PaHq} provides us with a $G$-equivariant biholomorphic map $\varphi: U \to Z$ to a quasi-projective algebraic $G$-variety with good quotient $\pi_Z: Z \to U\hq G$. Note that the proof of \cite[Thm.\ 3]{PaHq} heavily uses the assumption on the projectivity (and not just completeness) of $U\hq G$; see also the detailed discussion in Section~\ref{sect:embeddingGspaces} of this paper. The map $\varphi$ is induced from a holomorphic section in $(\pi_*\mathscr{F}^h)^G$ for some coherent algebraic $G$-sheaf $\mathscr{F}$ on $U$. Here, $\cdot^h$ denotes analytification of algebraic sheaves. Using the Coherence Theorem, an \'etale slice theorem (proven in Section~\ref{subsect:etaleslice}), and GAGA properties of the quotient $U\hq G$ we show that generically on $U\hq G$ the sheaf $(\pi_*\mathscr{F}^h)^G$ is algebraically isomorphic to the corresponding sheaf $(\pi_*\mathscr{F}^G)^h$ of equivariant \emph{algebraic} maps from $U$ to $Z$. Algebraicity of $\varphi$ follows. Consequently, $\pi$ is just a model of $\pi_Z$, hence itself a good quotient.

Several examples showing the necessity of the completeness and projectivity assumptions on the quotient $U\hq G$ are given in the individual sections (see Examples~\ref{ex:strangeopenset} and \ref{ex:FatouBieberbach}, Example~\ref{ex:nonalgebraicautomo}, Remark~\ref{rem:mumfordexample}, Remark~\ref{rems:finalremarks}.5, and Example~\ref{ex:nonprojective}).

With the help of our main result we approach a second important question in GIT: how many different open subsets $U$ with good quotient or analytic Hilbert quotient exist in a given $G$-variety? Using Theorem~\ref{thm:mainthmHamiltonian} as well as a finiteness result of Bia{\l}ynicki-Birula we prove
\begin{cor11.4}[Finiteness of momentum map quotients]
Let $G=K^\C$ be a complex reductive Lie group and let $X$ be a $G$-irreducible algebraic $G$-variety with at worst $1$-rational singularities. Then there exist only finitely many (necessarily Zariski-open) subsets of $X$ that can be realised as the set of $\mu$-semistable points with respect to some $K$-invariant K\"ahler structure and some momentum map $\mu: X \to \Lie(K)^*$ with compact zero fibre $\mu^{-1}(0)$.
\end{cor11.4}
Finally, we refine the statement of Theorem~\ref{thm:mainthmHamiltonian} in the case of semisimple group actions on projective varieties: using arguments of Hilbert-Mumford type and an ampleness criterion of Hausen and Bia{\l}ynicki-Birula/{\'S}wi{\polhk{e}}cicka we show the following result:
\begin{cor11.6}
Let $G$ be a connected semisimple complex Lie group with maximal compact subgroup $K$, and let $X$ be an irreducible projective algebraic Hamiltonian $G$-variety with at worst $1$-rational singularities with momentum map $\mu: X \to \Lie (K)^*$. Then there exists an ample $G$-linearised line bundle $L$ on $X$ such that $X(\mu)$ coincides with the set $X(L, G)$ of semistable points with respect to $L$.
\end{cor11.6}

In the final section we discuss how the results proven here fit into the general picture of Geometric Invariant Theory, momentum geometry, and their usage in the study of reductive group actions.

In an appendix we prove a technical result about the Luna stratification used in the proof of Theorem~\ref{thm:algebraicmap}.

\part{Tools}
\section{Notation and preliminaries on analytic Hilbert quotients}\label{sect:notations}
\subsection{Complex spaces and algebraic varieties}
In the following, a \emph{complex space} refers to a not necessarily reduced complex space with countable topology.
For a given complex space $X$ the structure sheaf is denoted by $\mathscr{H}_X$. Analytic
subsets are assumed to be closed. If $Z \hookrightarrow X$ is a closed complex subspace, its support is denoted by $|Z|$. By an \emph{algebraic variety} we mean an algebraic variety
defined over the field $\C$ of complex numbers. The structure sheaf of an algebraic variety
$X$ is denoted by $\mathscr{O}_X$. The complex space associated with a given algebraic variety $X$ is denoted by $X^h$. Given a morphism $\phi: X \to Y$ of algebraic varieties, we sometimes write $\phi^h: X^h \to Y^h$ for the induced map of complex spaces. If $\mathscr{F}$ is a sheaf of $\mathscr{O}_X$-modules on $X$, the associated sheaf of $\mathscr{H}_X$-modules will be denoted by $\mathscr{F}^h$. An \emph{open} subset of an algebraic variety always refers to a set that is open in the topology of $X^h$ and a \emph{Zariski-open} subset means a subset that is open in the algebraic Zariski-topology of $X$.
\subsection{Actions of Lie groups and analytic Hilbert quotients}\label{subsect:Actionsandquotients}
If $G$ is a real Lie group, then a \emph{complex $G$-space $Z$} is a
complex space with a real-analytic action $\alpha: G \times Z \to Z$ such that $G$ acts continuously on the structure sheaf $\mathscr{H}_X$ with its canonical Fréchet topology. If $G$ is a complex Lie group, a
\emph{holomorphic $G$-space
$Z$} is a complex $G$-space such that $G$ acts holomorphically on $\mathscr{H}_X$. A complex $G$-space $X$ is called \emph{$G$-irreducible} if $G$ acts transitively on the set of irreducible components of $X$. For more information on these notions, see \cite{HeinznerCoherent}.
Let $G$ be a complex reductive Lie group and $X$ a holomorphic $G$-space. A
complex space $Y$ together with a $G$-invariant surjective holomorphic map $\pi: X \to Y$ is called
an \emph{analytic Hilbert quotient} of $X$ by the action of $G$ if
\vspace{-2.5mm}
\begin{enumerate}
 \item $\pi$ is a locally Stein map, and
 \item $(\pi_*\mathscr{H}_X)^G = \mathscr{H}_Y$ holds.
\end{enumerate}
\vspace{-1mm}
Here, \emph{locally Stein} means that there exists an open covering of $Y$ by open Stein subspaces
$U_\alpha$ such that $\pi^{-1}(U_\alpha)$ is a Stein subspace of $X$ for all $\alpha$; by
$(\pi_*\mathscr{H}_X)^G$ we denote the sheaf $U \mapsto \mathscr{H}_X(\pi^{-1}(U))^G = \{f \in
\mathscr{H}_X(\pi^{-1}(U)) \mid f \;\; \text{is } G\text{-invariant}\}$, $U$ open in
$Y$.

An analytic Hilbert quotient of a holomorphic $G$-space $X$ is unique up to biholomorphism once it exists, and we will denote it by $X\hq G$. It has the following properties (see \cite{SemistableQuotients} and \cite{HeinznerCoherent}):

\begin{enumerate}
\item Given a $G$-invariant holomorphic map $\phi: X \to Z$ to a
complex space $Z$, there exists a unique holomorphic map $\widebar \phi: X\hq G \to Z$ such that
$\phi = \widebar \phi \circ \pi$.

\item For every Stein subspace $A$ of $X\hq G$ the inverse image $\pi^{-1}(A)$ is a Stein subspace
of $X$.
\item If $A_1$ and $A_2$ are $G$-invariant analytic (in particular, closed) subsets of $X$, then we have
$\pi(A_1) \cap \pi(A_2) = \pi(A_1 \cap A_2)$.
\item For a $G$-invariant closed complex subspace $A$ of $X$, which is defined by a $G$-invariant sheaf $\mathscr{I}_A$ of ideals, the image sheaf $(\pi_*\mathscr{I}_A)^G$ endows the image $\pi(A)$ in
$X\hq G$ with the structure of a closed complex subspace of $X\hq G$. Moreover, the restriction of $\pi$ to $A$ is an analytic Hilbert quotient for the action of $G$ on $A$.
\end{enumerate}
It follows that two points $x,x' \in X$ have the same image in $X\hq G$ if and only if
$\overline{G\acts x} \cap \overline{G\acts x'} \neq \emptyset$. For each $q \in X\hq
G$, the fibre $\pi^{-1}(q)$ contains a unique closed $G$-orbit $G\acts x$. The orbit $G\acts x$ is affine (see \cite[Prop. 2.3 and 2.5]{Snow}) and hence, the stabiliser $G_x$ of $x$ in $G$ is a complex reductive Lie group by a result of Matsushima \cite{Matsushima}.

If $X$ is a Stein space, then the analytic Hilbert quotient exists and has the properties listed above (see \cite{HeinznerGIT} and \cite{Snow}).

Motivated by the above, we introduce the following notation: if $X$ is a holomorphic $G$-space and $A$ is a $G$-stable subset of $X$, then we set \[\mathcal{S}_G^X(A) \definiere \{ x \in X \mid \overline{G \acts x} \cap A \neq \emptyset\}.\]We call $\mathcal{S}_G^X(A)$ the \emph{saturation of $A$ in $X$ (with respect to the $G$-action)}. If the context is clear, we will sometimes omit the superscript. If the analytic Hilbert quotient $\pi: X \to X\hq G$ exists and if $A$ is a $G$-invariant analytic subset of $X$, the results summarised in the previous paragraphs imply that $\mathcal{S}_G^X(A)  = \pi^{-1}(\pi(A))$.

\subsection{Good quotients}
Let $G$ be a complex reductive group endowed with its natural linear-algebraic structure. An
\emph{algebraic $G$-variety} is an algebraic variety $X$ together with an action of $G$ on $X$
such that the action map $G \times X \to X$ is regular. An algebraic $G$-variety is called \emph{$G$-irreducible} if $G$ acts transitively on the set of irreducible components of $X$.

An algebraic variety
$Y$ together with a $G$-invariant surjective regular map $\pi: X \to Y$ is called an \emph{algebraic
Hilbert quotient}, or \emph{good quotient}, of $X$ by the action of $G$  if
\vspace{-2.5mm}
\begin{enumerate}
 \item $\pi$ is affine, and
  \item $(\pi_*\mathscr{O}_X)^G = \mathscr{O}_Y$ holds.
\end{enumerate}

If $X$ is an algebraic $G$-variety, the associated complex space $X^h$ is in a natural way a holomorphic $G$-space. If the algebraic Hilbert quotient $\pi: X \to X\hq G$ exists, the associated holomorphic map $\pi^h: X^h \to (X\hq G)^h$ is an analytic Hilbert quotient for the action of $G$ on $X^h$, cf.\ \cite{Lunaalgebraicanalytic}.
\section{Technical tools}
\subsection{Chow and Rosenlicht quotients}\label{subsect:Rosenlicht}
Given an algebraic group $G$ and an algebraic $G$-variety $X$, there exists a natural construction which yields a geometric quotient $U_R/G$ of a $G$-invariant Zariski-open dense subset $U_R$ of $X$. One of the main ideas of our study is to compare the analytic Hilbert and momentum map quotients under discussion with this geometric quotient. Here we recall the relevant details.

A \emph{geometric quotient} for the action of an algebraic group $G$ on a variety $X$ is a surjective regular map $p: X \rightarrow Y$ having the following properties:
\begin{enumerate}
 \item For all $x \in X$, we have $p^{-1} (p(x)) = G \acts x$,
 \item $Y$ has the quotient topology with respect to $p$,
 \item $(\pi_*\mathscr{O}_X)^G = \mathscr{O}_Y$.
\end{enumerate}
For an irreducible projective algebraic variety $X$ and for $m, d \in \N$ let $\mathscr{C}_{m,d}(X)$ be the Chow variety of cycles of dimension $m$ and degree $d$ in $X$.
One of our main technical tools is the following result, which has appeared in a number of instances in the literature, e.g.\ see \cite{KapranovChow}, \cite{Lieberman}, \cite[Chap.\ 13]{BBSurvey} and references therein.
\begin{prop}\label{Chow}
Let $G$ be a connected algebraic group and let $X$ be an irreducible projective algebraic $G$-variety. Then there exist natural numbers $m,d
\in \N$, a $G$-invariant rational map $\varphi: X \dasharrow \mathscr{C}_{m,d}(X)$, and a
$G$-invariant Zariski-open subset $U_R \subset \dom(\varphi)$, called a \emph{Rosenlicht set} of $X$, such that
\begin{enumerate}
\item $U_R \subset X_{\mathrm{gen}} := \{x \in X \mid \dim G\acts x = m \text{ maximal}\}$
\item for all $u \in U_R$, we have $\varphi (u) = \overline{G\acts u}$ , considered as a (reduced) cycle of $X$,
\item the restriction $\varphi|_{U_R}: U_R \to \varphi(U_R)$ is a geometric quotient, $U_R/G \definiere \varphi(U_R)$ is normal and $\C(U_R/G) = \C(X)^G$.
\end{enumerate}
We call $\mathscr{C}_G(X) \definiere \varphi (U_R) = U_R/G$ a \emph{Chow quotient} or \emph{Rosenlicht quotient} of $X$ by $G$.
\end{prop}
\begin{figure}[h]
\begin{center}
  \includegraphics[height=4cm]{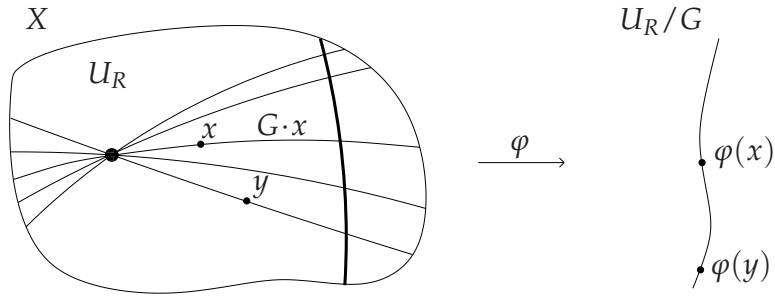}
  \caption{The Rosenlicht quotient}
\end{center}
\end{figure}
In Figure~1 the Rosenlicht set $U_R$ is the complement of the fat point and the thick line. Note that for two distinct points $p_1, p_2 \in U_R/G$ and for the corrresponding disjoint orbits $G\acts x_i = \varphi^{-1}(x_i) \subset U_R$, $i=1,2$, we might have $\overline{G\acts x_1}^X \cap \overline{G\acts x_2}^X \neq \emptyset$.

In greater generality, it has been shown by Rosenlicht~\cite{Rosenlicht2}, that every irreducible (not necessarily projective) algebraic $G$-variety contains a Zariski-open and dense subset $U_R$ such that a geometric quotient $U_R \to U_R/G$ with $\C(U_R/G) = \C(X)^G$ exists. Also, in this general case we call $U_R$ a \emph{Rosenlicht set} of $X$.
\subsection{Sumihiro neighbourhoods}\label{subsect:Sumihiro}
For local considerations it is often convenient to work with $G$-invariant subsets of projective spaces where $G$ acts linearly. This is the motivation for the following definition: Let $G$ be a connected algebraic group and let $X$ be an algebraic
$G$-variety. A \emph{Sumihiro neighbourhood} of a point $x \in X$ is a
$G$-invariant, Zariski-open, quasi-projective neighbourhood of $x$ in
$X$ that can be $G$-equivariantly embedded as a Zariski locally closed subset
of the projective space $\P (V)$ associated with some rational $G$-representation $\rho: G \to \mathrm{GL}(V)$.

In a normal irreducible algebraic $G$-variety, every point has a Sumihiro neighbourhood by a result of Sumihiro \cite{completion}.
\subsection{The Luna stratification}\label{subsect:Luna}
Let $G$ be a complex reductive Lie group and let $X$ be a holomorphic $G$-space. Let $G_x$, $G_y$ be two
isotropy subgroups of the action of $G$ on $X$. We define a preorder on the set of $G$-orbits in $X$
as follows: we say $G\acts x \leq G \acts y$ if and only if there exists an element $g \in G$ such that
$\Int(g)G_y < G_x$. Here, $\Int (g): G \rightarrow G$ is given by $\Int (g)(h) = ghg^{-1}$. This induces an equivalence relation on the set of $G$-orbits: $G\acts x \sim G \acts y$
if and only if there exists an element $g \in G$ such that $\Int (g)G_y = G_x$. We denote the equivalence class of an orbit $G\acts x$ by $\Type(G\acts x)$ and call it the \emph{orbit type} of $x$.

Assume now that the $G$-action on $X$ admits an analytic Hilbert quotient $\pi: X \to X\hq G$. With the help of the above, we define an equivalence relation on the points of $X\hq G$: two points $p, q \in X\hq G$
are defined to be equivalent if the uniquely defined closed orbits $G\acts x$ and $G\acts y$ in $\pi^{-1}(p)$ and $\pi^{-1}(q)$, respectively, have the same orbit type. If $H$ is a representative of a conjugacy class of reductive subgroups of $G$, we denote the corresponding equivalence class by $(X\hq G)^{(H)}$.

The following proposition generalises a result of Luna \cite{LunaSlice} from the algebraic to the analytic category. The smooth analytic case is also handled in \cite{SjamaarSlices}.
\begin{prop}[\cite{Extensionofsymplectic}]
Let $X$ be a $G$-irreducible complex $G$-space with analytic Hilbert quotient $\pi: X \to X\hq G$. Then the following holds:
\begin{enumerate}
\item There exists a uniquely determined equivalence class $(X\hq G)^{\mathrm{princ}}$, called the \emph{principal stratum}, corresponding to the minimal conjugacy class $(H)$ of isotropy groups of closed orbits in $X$, which is analytically Zariski-open and dense in $X\hq G$.
\item If we set $X^{\mathrm{princ}} \definiere \pi^{-1}((X\hq G)^{\mathrm{princ}})$, then the restriction of $\pi$ to $X^{\mathrm{princ}}$, \[\pi|_{X^{\mathrm{princ}}}: X^{\mathrm{princ}} \to (X\hq G)^{\mathrm{princ}},\]
    is a holomorphic fibre bundle with typical fibre $G \times _H Z$, where $Z$ is an affine algebraic $H$-variety with $\C[Z]^H = \C$.
\end{enumerate}
\end{prop}
Here, $G \times_H Z$ denotes the quotient of $G \times Z$ by the diagonal $H$-action given by $h \acts (g,z) = (gh^{-1},h\cdot z)$.
\subsection{Reduction to the case of connected groups}
The following considerations will allow us to restrict our attention to action of connected groups in the proof of Theorem~\ref{thm:mainthm}.
\begin{lemma}\label{lem:connectedisenough}
Assume that the statement of Theorem~\ref{thm:mainthm} has been shown for all \emph{connected} complex reductive groups $G$. Then Theorem~\ref{thm:mainthm} holds for every complex reductive group $G$.
\end{lemma}
\begin{proof}
Since $G_0$ is a reductive subgroup of $G$, the analytic Hilbert quotient $\pi_{G_0}: U \to U\hq G_0$ exists, cf.\ \cite[2.\ Cor.]{SemistableQuotients}. Furthermore, we obtain the following commutative diagram:
\[\begin{xymatrix}{
U \ar[r]^\pi\ar[d]_{\pi_{G_0}}& U\hq G\\
U\hq G_0\ar[ru]_{\pi_C}.&
}
\end{xymatrix}\]
Here, $\pi_C$ is the algebraic Hilbert quotient for the induced action of the finite group $C
\definiere G/G_0$ on $U \hq G_0$. Consequently, $U\hq G_0$ is projective. Since $X$ is assumed to be normal, the irreducible components of $U$ are disjoint. The connected component $G_0$ of the identity of $G$ stabilises each irreducible component. Observe that the restriction of $\pi_{G_0}$ to any of the irreducible components $U_j = U \cap X_j$ of $U$, $\pi_j \definiere \pi_{G_0}|_{U_j}: U_j \to \pi(U_j) \cong U_j \hq G$, is an analytic Hilbert quotient. As an analytic subset of $U\hq G_0$, the quotient $U_j\hq G_0$ is likewise a projective algebraic variety. Furthermore, note that $\pi_i(U_i) \cap \pi_j(U_j) = \emptyset$ for $i\neq j$. By assumption, $U_j$ is Zariski-open and dense in $X_j$. It follows that $U$ is Zariski-open and dense in $X$. Again by assumption, each $\pi_j$ is a good quotient. It follows that $\pi_{G_0}$ is a good quotient. It is a classical result that the quotient of a projective variety by a finite group is a good quotient. As a consequence $\pi$ is an algebraic Hilbert quotient.
\end{proof}
\section{Embedding holomorphic $G$-spaces}\label{sect:embeddingGspaces}
In this section we shortly discuss the second main result of \cite{PaHq} as well as those technical parts of the proof that will be important for our study here.

The following "Algebraicity Theorem" states that on a given complex $G$-space with projective algebraic quotient there exists an essentially unique algebraic structure making the group action algebraic and the quotient a good quotient in the sense of GIT.
\begin{thm*}[Algebraicity Theorem, Thm.~3 of \cite{PaHq}]
Let $G$ be a complex reductive Lie group. Let $X$ be a holomorphic $G$-space such that the analytic
Hilbert quotient $\pi: X \rightarrow X\hq
G$ exists and such that $X\hq G$ is projective algebraic. Then, up to $G$-equivariant algebraic isomorphisms, there exists a uniquely determined
quasi-projective algebraic $G$-variety $Z$ with algebraic Hilbert quotient $\pi_Z: Z \rightarrow
Z\hq G$ such that $X$ is $G$-equivariantly biholomorphic to $Z$.
\end{thm*}
Let $V$ be a finite-dimensional $G$-module, and let $\mathscr{L}$ be the locally free sheaf associated with an ample line bundle $L$ on $X\hq G$. Then, for any $m \in \N$, the sheaf $V \otimes \pi^*\mathscr{L}^{\otimes m}$ is the sheaf of sections of a holomorphic $G$-vector bundle over $X$. Invariant sections of this bundle can be used to construct linearly equivariant maps on $X$ by the following
\begin{lemma}[\cite{PaHq}, Lem.\ 8.1]\label{lem:basicconstruction}
Let $G$ be a complex reductive Lie group and let $X$ be a complex $G$-space such that the analytic Hilbert quotient $X\hq G$ exists as a projective algebraic complex space. Let $V$ be a $G$-module, let $L$ be an algebraic line bundle on $X\hq G$ and let $\mathscr{L}$ be the associated locally free sheaf.
Then, an element $s \in H^0\bigl(X\hq G, \pi_*(V\tensor \pi^* \mathscr{L})^G\bigr)$ yields
a $G$-equivariant holomorphic map $\sigma: X \rightarrow \mathcal{V}$
into the algebraic $G$-vector bundle $\mathcal{V}=  V \tensor L$ over $X\hq G$.
\end{lemma}
The proof of the Algebraicity Theorem in \cite{PaHq} shows that given an ample line bundle $L$ on $Q$ there exist finitely many $G$-modules $V_1, \dots, V_k$, natural numbers $m_1, \dots, m_k \in \N$, and $G$-invariant sections $s_j \in H^0\bigl(X\hq G, \pi_*(V_j \otimes \pi^*\mathscr{L}^{\otimes m_j})^G\bigr)$ such that the direct product $\Psi \definiere \bigoplus_{j=1}^k \sigma_j: X \to \mathcal{V} \definiere \bigoplus_{j=1}^k (V_j \otimes L^{\otimes m_j})$ of the maps provided by Lemma~\ref{lem:basicconstruction} yields a $G$-equivariant holomorphic embedding of $X$ into $\mathcal{V}$ (Embedding Theorem, Thm.\ 8.6) with Zariski-closed image $ \Psi (X)$ (Thm.\ 9.3). Furthermore, the good quotient of the natural $G$-action on $\Psi(X)$ exists and $\Psi(X)\hq G$ is naturally biregular to $X\hq G$.
\begin{rem}
In contrast to many of the results proven in this paper, the Algebraicity Theorem discussed above heavily relies on the assumption on the projectivity (and not just completenes) of the quotient $Q$. The assumption is mainly used to construct sections of $\pi_*(V_j \otimes \pi^* \mathscr{L}^{\otimes m_j})^G$ using Serre vanishing for the ample bundle $L$.
\end{rem}
\section{Meromorphic functions on analytic Hilbert quotients}\label{sect:meromorphic}
In this section we study the function fields of analytic Hilbert quotients. In particular, we study conditions under which invariant meromorphic functions on $X$ descent to the quotient $X\hq G$.
\subsection{Meromorphic functions and meromorphic graphs}
We denote the sheaf of germs of \emph{meromorphic functions} on a not necessarily reduced complex space $X$ by $\mathscr{M}_X$. Let $f \in \mathscr{M}_X(X)$ be a meromorphic function on $X$. The \emph{sheaf of denominators} of $f$ is the sheaf $\mathscr{D}(f)$ with stalks
$\mathscr{D}(f)_x = \{v_x \in \mathscr{H}_x : v_x f_x \in \mathscr{H}_x \}$. We
define the \emph{polar variety} of $f$ to be the closed complex subspace defined by $\mathscr{D}(f)$ and denote it by $P_f$. The polar variety is the smallest subset of $X$ such that $f$ is holomorphic on $X \setminus P_f$.

We will describe meromorphic functions via their graphs, see Proposition~\ref{prop:graphtofunction} below. Consider a closed complex subspace $\Gamma \hookrightarrow X \times \P_1$ and denote the canonical projection to $X$ by $\sigma \definiere \textrm{pr}_1|_\Gamma: \Gamma \to X$. Then $\Gamma$ is called a \emph{holomorphic graph at }$p\in X$, if there exists an open neighbourhood $U$ of $p$ in $X$ such that
\begin{enumerate}
\item $\sigma|_{\sigma^{-1}(U)}: \sigma^{-1}(U) \to U$ is biholomorphic,
\item $\sigma^{-1}(U) \cap (U \times \{\infty\}) = \emptyset$.
\end{enumerate}
Clearly, the graph of a holomorphic function $f \in \mathscr{H}_X(X)$ is a holomorphic graph at every point $p \in X$. A closed complex subspace $\Gamma \hookrightarrow X \times \P_1$ with canonical map $\sigma: \Gamma \to X$ is called a \emph{meromorphic graph over }$X$, if there exists an analytic set $A \subset X$ with the following properties:
\begin{enumerate}
\item $A$ and $\sigma^{-1}(A)$ are analytically rare,
\item $\Gamma$ is a holomorphic graph outside $A$.
\end{enumerate}
Recall that an analytic subset $Z$ in a complex space $X$ is called \emph{analytically rare} if for every open subset $U$ of $X$ the restriction map $\mathscr{H}_X(U) \to \mathscr{H}_X(U \setminus Z)$ is injective. Please consult \cite[0.43]{Fischer} for further details on this notion.
\begin{prop}[\cite{FischerMeromorphic}]\label{prop:graphtofunction}
The graph $\Gamma_f$ of a meromorphic function on a complex
space is a meromorphic graph. Conversely, let $\Gamma \hookrightarrow X \times \P_1$ be a meromorphic graph over $X$. Then there exists a uniquely determined meromorphic function $f \in \mathscr{M}_X(X)$ such
that $\Gamma = \Gamma_f$.
\end{prop}
\subsection{Descent of invariant meromorphic functions}
The following is one of the main technical tools used in the sequel.
\begin{prop}\label{prop:meromorphicdescent}
Let $G$ be a complex reductive Lie group and let $X$ be a holomorphic $G$-space with analytic Hilbert quotient $\pi: X \to Q \definiere X\hq G$. Let $P \hookrightarrow X$ be a $G$-invariant $\pi$-saturated closed subspace such that $|P|$ is analytically rare in $X$. Let $f \in \mathscr{H}_X(X \setminus |P|)^G$ and assume that $f$ extends to a $G$-invariant meromorphic function $F \in \mathscr{M}_X(X)^G$. Then the uniquely determined holomorphic function $\hat f \in \mathscr{H}_{X/\negthickspace/ G}(Q \setminus \pi(|P|))$ that fulfills $\pi^*(\hat f) = f$ extends to a meromorphic function $\hat F \in \mathscr{M}_Q(Q)$ with $\pi^*(\hat F) =F$.
\end{prop}
\begin{proof}
The map $\Pi \definiere \pi \times \id_{\P_1}: X \times \P_1 \to Q \times \P$ is an analytic Hilbert quotient for the action of $G$ on $X \times \P_1$ induced by its action on the first factor. Let $F \in \mathscr{M}_X(X)^G$ be as in the claim and let $\Gamma \hookrightarrow X\times \P_1$ be its graph. Endow $\widehat \Gamma \definiere \Pi(\Gamma) \subset Q \times \P_1$ with the canonical complex structure, cf.\ Section~\ref{subsect:Actionsandquotients}, such that $\Phi \definiere \Pi|_\Gamma : \Gamma \to \widehat \Gamma$ is an analytic Hilbert quotient.This leads to the following commutative diagram:
\begin{equation*}
\begin{xymatrix}{
 X \times \mathbb{P}_1\ar[d]^{\Pi} &\ar@{_{(}->}[l]
\Gamma\ar[r]^{\sigma_\Gamma}\ar[d]^{\Phi}  & X\ar^{\pi}[d] \\
Q \times \mathbb{P}_1 & \ar@{_{(}->}[l] \widehat \Gamma_{\ } \ar[r]^{\sigma_{\widehat \Gamma}}  & Q.
}
\end{xymatrix}
\end{equation*}
The analytic subset $\pi(|P|) \subset Q$ is analytically rare. Furthermore, by assumption, $\widehat \Gamma$ is a holomorphic graph over $Q \setminus \pi(|P|)$. To establish the claim it suffices to show that the set $\sigma_{\widehat \Gamma}^{-1}(\pi(|P|))$ is analytically rare in $\widehat \Gamma$. This follows from $\Phi^{-1}(\sigma_{\widehat \Gamma}(\pi(|P|))) = \sigma_\Gamma^{-1}(|P|)$ and from the fact that the latter set is analytically rare in $\Gamma$ by assumption.
\end{proof}
\begin{cor}\label{cor:meromorphicdescentgenericcase}
Let $X$ be a reduced $G$-irreducible holomorphic $G$-space with analytic Hilbert quotient $\pi: X\to Q \definiere X\hq G$. Assume that there exists a point $x_0 \in X_{\mathrm{gen}}$ such that $G\acts x_0$ is closed in $X$. Then for every $f \in \mathscr{M}_X(X)^G$ there exists a $\bar f \in \mathscr{M}_Q(Q)$ such that $f = \pi^* (\bar f)$.
\end{cor}
\begin{proof}
We have to find $P$ as in the statement of Proposition~\ref{prop:meromorphicdescent}. Set $E \definiere \pi^{-1}(\pi (X \setminus X_{\mathrm{gen}}))$ and $X^s \definiere X \setminus E$. It follows from the existence of $x_0$ that $E$ is a nowhere dense analytic subset of $X$. Note that $X^s$ is $\pi$-saturated in $X$ and that the restriction of $\pi$ to $X^s$ is a geometric quotient for the $G$-action on $X^s$, cf.\ \cite[Lem.\ 3.4]{PaHq}. It follows that $P \definiere E \cup \pi^{-1}(\pi (P_f))$ is a nowhere dense saturated analytic subset of $X$ and hence analytically rare in $X$. Furthermore, $f|_{X\setminus P} \in \mathscr{H}_X (X\setminus P)$ by construction. Hence, the claim follows from Proposition~\ref{prop:meromorphicdescent}.
\end{proof}
\begin{rem}
Special cases of Proposition~\ref{prop:meromorphicdescent} and Corollary~\ref{cor:meromorphicdescentgenericcase} have been proven in \cite[Sect.\ 7.2]{PaHq} for reduced Hamiltonian $G$-varieties and for quotients arising from semi\-sta\-ble points of momentum maps.
\end{rem}
\part{Proof of the main result}
In this part, we prove Theorem~\ref{thm:mainthm}. We have separated three preliminary steps from the main argument. These are given in Sections~\ref{sect:Zariskiopen}, \ref{sect:algebraicmap}, and \ref{sect:coherence}, respectively. The proof of Theorem~\ref{thm:mainthm} is then given in Section~\ref{sect:maintheoremproof}.

Due to Lemma \ref{lem:connectedisenough} we may assume that $G$ is connected throughout this part. Furthermore, we assume that $U$ is non-empty. For an analytic subset $A$ of $X$ we set $A^{ss} \definiere A\cap U$.
\section{The set $U$ is open in the Zariski-topology}\label{sect:Zariskiopen}
In this section, we will prove the following result:
\begin{thm}[Openness Theorem]\label{thm:Zariskiopen}
Let $G$ be a connected complex reductive Lie group and let $X$ be an irreducible normal algebraic $G$-variety. Let $U$ be a $G$-invariant analytically Zariski-open subset of $X$ such that the analytic Hilbert quotient $\pi: U \to U \hq G$ exists. If $U\hq G$ is a complete algebraic variety, then $U$ is Zariski-open in $X$.
\end{thm}
\subsection{Generic closed orbits}\label{subsect:genericorbits}
Recall that for an irreducible algebraic $G$-variety $X$, we have defined
$ X_{\mathrm{gen}} = \{x \in X \mid \dim G\acts x = m\}$,
where $m= \max_{x \in X}\{\dim G\acts x\}$. The behaviour of a quotient $\pi: X \to X\hq G$ is particularly easy to control if the generic $G$-orbit is closed in $X$. The following result allows to reduce some considerations to this situation.
\begin{lemma}\label{lem:reductiontogeneric}
Let $G$ be a connected complex reductive Lie group. Let $X$ be an irreducible algebraic $G$-variety and $U$ a $G$-invariant analytically Zariski-open subset of $X$ such that the analytic Hilbert quotient $\pi: U \to Q $ exists. Then there exists a $G$-stable irreducible subvariety $Y$ of $X$ such that $\pi(Y^{ss}) = Q$ and such that $Y_{\mathrm{gen}}\cap U$ contains an orbit that is closed in $U$.
\end{lemma}
\begin{proof}
We argue by induction on the dimension of $X$. If $\dim X = 0$, there is nothing to prove. In the  general case, if $X_{\mathrm{gen}}$ contains an orbit that is closed in $U$, we are done. So we can assume that all closed orbits of $U$ lie in $Z \definiere X \setminus X_{\mathrm{gen}}$. It follows that $\pi(Z \cap U)=Q$. Let $Z = \bigcup_{j = 1}^m Z_j$ be the decomposition of $Z$ into irreducible components. Then there exists $j_0 \in \{1, \dots, m\}$ such that $Q = \pi(Z_{j_0}\cap U) \cong (Z_{j_0}\cap U)\hq G$. Here, $Z_{j_0}\cap U$ is analytically Zariski-open in $Z_{j_0}$. Since $\dim Z_{j_0} < \dim X$, induction applies.
\end{proof}
\subsection{Proof of Theorem~\ref{thm:Zariskiopen}}\label{subsect:proofofZariskiopen}
We start by studying the case where the generic $G$-orbit is closed in $U$.
\begin{prop}\label{prop:rosenlichtandsemistable}
Let $G$ be a connected complex reductive Lie group, let $X$ be an irreducible algebraic $G$-variety, and let $U \subset X$ be a $G$-invariant analytically Zariski-open subset of $X$ such that the analytic Hilbert quotient $\pi: U \to Q$ exists. Assume that $Q$ is a complete algebraic variety and that $U \cap X_{\mathrm{gen}}$ contains an orbit that is closed in $U$.

If $V_R$ is any Rosenlicht subset of $X$, then $U\cap V_R$ contains a $G$-invariant Zariski-open (Rosenlicht) subset $U_R$ of $V_R$ consisting of $G$-orbits that are closed in $U$ such that there exists an open algebraic embedding $\imath: U_R/G \hookrightarrow Q$ making the following diagram commutative:
\[\begin{xymatrix}{
  &\ar_{\pi}[ld]  U_R  \ar^{p}[dr]&  \\
Q &     & \ar@{_{(}->}_{\imath}[ll] U_R/G.}
\end{xymatrix}
\]
In particular, the quotient map $\pi: U \to Q$ extends to a rational map $X \dasharrow Q$.
\end{prop}
\begin{proof}
Let $V_R$ be any Rosenlicht set of $X$ and let $p: V_R \to V_R/G$ be the quotient map. Without loss of generality we can assume that $V_R/G$ is affine. Our first aim is to show that $V_R/G$ is birational to $Q$. By Corollary~\ref{cor:meromorphicdescentgenericcase}, the restriction $f|_U$ of every $G$-invariant meromorphic function $f\in \mathscr{M}_X(X)^G$ to $U$ descends to a meromorphic function on $Q$. As a consequence of Grothendieck's generalisation of Serre's GAGA results, cf.\ \cite[XII.4.]{SGA1}, these meromorphic functions are in fact rational; see also \cite[Thm.\ VIII.3.1.1]{ShafarevichII}. Consider the rational map $\varphi: Q \dasharrow V_R/G$ that corresponds to $p|_U$. Since elements of $\varphi^*(\C(V_R/G)) = \C(X)^G \subset \C(Q)$ separate orbits in $V_R$, they separate orbits in
\[U_R \definiere U \setminus \mathcal{S}^U_G(V^c_R \cap U) \subset U\cap V_R.\]

This set is non-empty, analytically Zariski-open, and $\pi$-saturated in $U$. Since $U$ contains a closed orbit of generic orbit dimension and since $Q$ is complete, $U_R$ is mapped to a non-empty Zariski-open subset of $Q$, as a consequence of GAGA, see \cite[XII.4.]{SGA1}. Therefore, $\varphi$ is generically one-to-one, hence birational. Let $\imath \definiere \varphi^{-1} : V_R/G \dasharrow Q$. Without loss of generality, $V_R/G$ coincides with the set where $\imath$ is an isomorphism onto its image. It follows that $\imath: V_R/G \hookrightarrow Q$ is an open embedding. We obtain the following commutative diagram
\[\begin{xymatrix}{
U_R  \ar@{^{(}->}[rr]\ar[d]_{p|_{U_R}}&  &U\ar[d]^\pi \\
p(U_R) \ar@{^{(}->}[rr]^<<<<<<<<<<{\imath|_{p(U_R)}}&  &Q.
}
  \end{xymatrix}
\]
Since $Q$ is complete, the image $p(U_R)$ is Zariski-open in $Q$. It follows that $U_R$ is Zariski-open in $V_R$ and contained in $U$. This shows the claim.
\end{proof}
We now return to the general case.
\begin{lemma}\label{lem:UcontainsZopen}
Let $G$ be a connected complex reductive Lie group, $X$ an irreducible $G$-variety, and $U$ a non-empty  $G$-invariant analytically Zariski-open subset of $X$ such that the analytic Hilbert quotient $\pi: U \to Q$ exist and $Q$ is a complete algebraic variety.
Assume that every point $x\in U$ whose orbit $G\acts x$ is closed in $U$ has a Sumihiro neighbourhood. Then, $U$ contains a non-empty $G$-invariant Zariski-open subset of $X$.
\end{lemma}
\begin{proof}
By Lemma~\ref{lem:reductiontogeneric}, there exists an irreducible $G$-invariant subvariety $Y$ of $X$ such that $Y^{ss}\hq G = Q$ and such that $Y^{ss} \cap Y_{\mathrm{gen}}$ contains an orbit that is closed in $U$. By Proposition~\ref{prop:rosenlichtandsemistable} there exists a $G$-invariant Zariski-open subset $A$ of $Y$ contained in $U$ such that $\pi(A)$ is an open subset of $X(\mu)\hq G$. Furthermore, $A$ consists of
orbits which are closed in $U$.

Let $W$ be an irreducible Sumihiro neighbourhood of a point $x \in A$ and let $\psi: W \to  \P(V)$ be a $G$-equivariant embedding of $W$ into the projective
space associated with a rational $G$-representation $V$. Let $Z$ be the closure of $\psi(W)$ in $\P(V)$.
Given a Rosenlicht subset $U_R$ of $Z$ as in Proposition~\ref{Chow}, Lemma 6.3 of \cite{PaHq} implies that $\mathcal{S}^Z_G(\psi(A \cap W))\cap U_R$ is constructible in
$U_R$. Therefore,
$\mathcal{S}_G^X(A \cap W)\cap W$ contains a $G$-invariant Zariski-open subset $\widetilde U$ of its closure. By construction, $\pi(A \cap W)$ is open in $Q$, and hence, $\pi^{-1}(\pi(A \cap W)) \cap W$ is an open subset of $X$ that is contained in
$\mathcal{S}_G(A \cap W)\cap W$. We conclude that $\widetilde{U}$ is Zariski-open in
$X$.
\end{proof}

\begin{prop}\label{prop:SumihiroZopen}
Let $G$ be a connected complex reductive Lie group, $X$ a $G$-variety, and $U$ a $G$-invariant analytically Zariski-open subset of $X$ such that the analytic Hilbert quotient $\pi: U \to Q$ exist and $Q$ is a complete algebraic variety.
Assume that every point $x \in U$ whose orbit $G\acts x$ is closed in $U$ has a Sumihiro neighbourhood. Then $U$ is Zariski-open in $X$.
\end{prop}
\begin{proof}
Let $X = \bigcup_j X_j$ be the decomposition of $X$ into irreducible components. Then, there exists a $j_0 \in \{1, \dots, m\}$ such that $X_{j_0}^{ss}$ is analytically Zariski-open in $X_{j_0}$ and non-empty. Let $V_{j_0}$ be the $G$-invariant Zariski-open subset of $X_{j_0}^{ss}$ whose existence is
guaranteed by the previous Lemma. Set $\widetilde {X} \definiere X \setminus \bigl(V_{j_0} \setminus \bigcup_{j\neq j_0} X_j \bigr)$. Then either $\widetilde{X} = U^c \subset X$ or $\widetilde X^{ss} \neq \emptyset$. In the first case, we are done, since $\widetilde {X}$ is algebraic in $X$. In the second case, we notice that $\widetilde X^{ss}$ is analytically Zariski-open and $G$-invariant in $\widetilde X$, that $\widetilde X^{ss} \hq G = \pi (\widetilde X^{ss}) \subset Q$ is a complete algebraic variety by GAGA \cite[Thm.\ XII.4.4]{SGA1}, and that the existence of Sumihiro neighbourhoods is inherited by $\widetilde X^{ss}$. We proceed by Noetherian induction.
\end{proof}
\begin{proof}[Proof of Theorem~\ref{thm:Zariskiopen}]
Using Proposition~\ref{prop:SumihiroZopen} it suffices to note that, since $X$ is assumed to be normal and $G$ to be connected, every point in $X$ has a Sumihiro neighbourhood, cf.\ Section~\ref{subsect:Sumihiro}.
\end{proof}
The following examples show that completeness of the quotient is a necessary assumption for Theorem~\ref{thm:algebraicmap}:
\begin{ex}\label{ex:strangeopenset}
Let $X= \C^2$, and let $\Phi$ be the non-algebraic, holomorphic automorphism of $\C^2$ given by $(z,w) \mapsto (w, e^z - w)$. Then, setting $U \definiere \C^2 \setminus \{e^w -z = 0 \}$ and $Y \definiere \C^2 \setminus\{ w=0 \}$, the map $\Phi$ induces a biholomorphic map $\varphi: U \to Y$. Note that $Y$ is an algebraic variety and that $\varphi$ is an analytic Hilbert quotient for the action of the trivial group on $U \subset \C^2$.
\end{ex}
In the previous example, the set $U$ was analytically, but not algebraically Zariski-open. The following example shows that, even worse, the complement of $U$ in $X$ can have non-empty interior if we drop the completeness assumption.
\begin{ex}\label{ex:FatouBieberbach}
Consider again $X\definiere \C^2$ and let $U \subset \C^2$ be a so-called \emph{Fatou-Bieberbach domain}, a proper (metrically open) subdomain of $\C^2$ admitting a biholomorphic map $\varphi: U \to \C^2$ (these domains arise in holomorphic dynamics; for concrete examples, see \cite[Ex.\ 6.3.2]{HolomorphicDynamics}). Then, we note that again $\varphi$ is an analytic Hilbert quotient for the action of the trivial group on $U \subset X$.
\end{ex}

\section{Algebraicity of the quotient map}\label{sect:algebraicmap}
The aim of this section is to prove the following result.
\begin{thm}[Algebraicity Theorem]\label{thm:algebraicmap}
Let $G$ be a connected complex reductive Lie group, let $X$ be an irreducible normal algebraic $G$-variety, and let $U$ be a $G$-invariant Zariski-open subset of $X$ such that the analytic Hilbert quotient $\pi: U \to Q$ exists. If $Q$ is a complete algebraic variety, then the quotient map $\pi$ is algebraic.
\end{thm}
\subsection{The principal Luna stratum}
Since the quotient map is a priori holomorphic, it sufffices to show that it extends to a rational map $\pi: X \dasharrow U\hq G$. Hence, it is enough to understand the behaviour of $\pi$ on an algebraically Zariski-open subset of $U$. We will use the principal stratum of the Luna stratification for this purpose.

Let $G,X$, and $\pi:U \to Q$ be as in the hypotheses of Theorem~\ref{thm:algebraicmap}. Let $S \definiere Q_{\mathrm{princ}}$ be the principal Luna stratum, let $Y$ be the closed $G$-invariant subvariety of $X$ introduced in Lemma~\ref{lem:reductiontogeneric}, and let $\pi_Y: Y^{ss} \to Q$ denote the restriction of $\pi$ to $Y^{ss}$.
\begin{lemma}\label{lem:principalorbitsclosed}
Under the assumptions listed above, $\pi_Y^{-1}(S)$ consists of closed orbits of $G$ in $Y^{ss}$.
\end{lemma}
\begin{proof}
It follows from the construction of $Y$ that the principal orbit type of $Q$ equals that of $Y^{ss}$. Let $H$ be a representative of the corresponding conjugacy class of reductive subgroups of $G$. As we have seen in Section~\ref{subsect:Luna}, the map $\pi_Y: \pi_Y^{-1}(S) \to S$ is a holomorphic fibre bundle with fibre $G/H$. Consequently, every orbit $G \acts y \subset \pi_Y^{-1}(S)$ is closed in $Y^{ss}$.
\end{proof}
The proof of the following result strongly resembles the proofs of the results in Section~\ref{sect:Zariskiopen}, and has therefore been deferred to Appendix~\ref{appendix:Luna}.
\begin{lemma}\label{lem:LunaZopen}
Under the assumptions listed above, $\pi^{-1}(S)$ is Zariski-open in $X$.
\end{lemma}
\subsection{Proof of Theorem \ref{thm:algebraicmap}}
In the course of the proof of Theorem~\ref{thm:algebraicmap} we will construct a "well-defined family of compact algebraic cycles" from the orbits of the group action on $X$. Since the definition of "algebraic family" is technically involved (Definition I.3.10 and Definition I.3.11 of \cite{KollarRatCurves}), we have chosen to state only the following simplified version of a criterion by Koll\'ar which will be used in our proof.
\begin{thm}[cf.\ Thm.~I.3.17 of \cite{KollarRatCurves}]\label{thm:algebraicfamily}
Let $X, W$ be normal varieties, let $U$ be an irreducible subvariety of $X \times W$, and let $\pi_2: U \to W$ be the natural projection.  Then, $(\pi_2: U \to W)$ is a well defined-family of $k$-dimensional compact algebraic cycles of $X$ over $W$ if $\pi_2$ is proper and surjective, and if every fibre of $\pi_2$ has dimension $k$.
\end{thm}
Given a well defined-family of $k$-dimensional compact algebraic cycles of a projective variety $X$ (with ample line bundle $\mathscr{O}_X(1)$) over an irreducible $W$, there exists a natural number $d$, and a uniquely determined morphism $\psi: W \to \mathscr{C}_{k,d}(X)$ into the Chow variety of cycles of dimension $k$ and degree $d$ in $X$ such that $U$ is isomorphic to the pullback (via $\psi$) of the universal family over $\mathscr{C}_{k,d}(X)$ , see \cite[Thm.~I.3.21]{KollarRatCurves}.

We will furthermore apply the following result of Hausen.
\begin{prop}[Prop.~2.6 of \cite{HausenGITwithWeildivisors}]\label{prop:qpinvariant}
Let $G$ be a connected linear algebraic group, let $X$ be a normal $G$-variety, and let $U\subset X$ be an open subset. If $U$ is quasi-projective, then $G\acts U$ is quasi-projective. In particular, the maximal quasi-projective open subset $U$ of $X$ is $G$-invariant.
\end{prop}
\begin{proof}[Proof of Theorem~\ref{thm:algebraicmap}]
Since we are only interested in the behaviour of $\pi$ on a Zariski-open subset of $X$, due to Proposition~\ref{prop:Lunaalgebraic}, we may assume in the following that $X=U=\pi^{-1}(S)$.

Consider the following algebraic map defined by the $G$-action on $X$:
\begin{align*}
\Psi: G \times X &\to X \times X \\
      (g, x) &\mapsto (x, g\acts x),
\end{align*}
and let $\mathcal{R} \definiere \Psi(G \times X)$ be the corresponding orbit relation in $X \times X$. Since $G$ is connected and $X$ is irreducible, its closure $\overline{\mathcal{R}}$ is a closed $(G\times G)$-invariant irreducible subvariety of $X\times X$. Let $\pi_1, \pi_2: X \times X \to X$ be the canonical projections to the first and second factor, respectively. Since $\pi_2$ is $G$-equivariant and surjective, there exists a non-empty $G$-invariant Zariski-open subset $\mathscr{U} \subset X$ such that
\begin{enumerate}
\item $\mathscr{U}$ is smooth,
\item $\dim (\pi_2|_{\overline{\mathcal{R}}})^{-1} (u) = \dim \overline{\mathcal{R}} - \dim X$,
\item $(\pi_2|_{\overline{\mathcal{R}}})^{-1} (u) = \overline{G\acts u} \times \{u\}$.
\end{enumerate}
We will show that there exists a Zariski-open $G$-invariant subset $\mathscr{U}'$ of $\mathscr{U}$ such that the restriction of $\pi$ to $\mathscr{U}'$ is a regular map.

First, we consider the intersection of fibres of $\pi_2|_{\overline{\mathcal{R}}}$ with the subvariety $Y$. By Lemma~\ref{lem:principalorbitsclosed}, the closure of every $G$-orbit in $X=\pi^{-1}(S)$ intersects $Y=Y^{ss}$ in a unique closed $G$-orbit $G\acts y$ with stabiliser $G_y$ that is conjugate to the minimal isotropy group $H$  in $G$. We set
\[\Gamma \definiere \overline{\mathcal{R}} \cap (Y \times \mathscr{U}).\]
This is a closed $G$-invariant subvariety of $Y \times \mathscr{U}$. Since the closure of every $G$-orbit in $\mathscr{U}$ intersects $Y$, we have $\pi_2(\Gamma) = \mathscr{U}$. There exists an irreducible component $\Gamma_0$ of $\Gamma$ such that $\pi_2(\Gamma_0)$ contains a $G$-invariant Zariski-open subset of $\mathscr{U}$. Hence, without loss of generality, we can assume that $\Gamma$ is irreducible and that $\pi_2(\Gamma)=\mathscr{U}$. Note that for all $u \in \mathscr{U}$ there exists some $y \in Y \cap \overline{G\acts u}$ such that $(\pi_2|_\Gamma)^{-1}(u) = G\acts y \times \{u\}$. As a next step we show how to build from $\Gamma$ an algebraic family of compact cycles in $Y$ with very controlled behaviour.

Let $Y'$ be the maximal quasi-projective subset of the smooth locus of $Y$. The subset $Y'$ is $G$-invariant by Proposition~\ref{prop:qpinvariant}. Note that the set $\Gamma' \definiere \Gamma \cap (Y' \times \mathscr{U})$ is non-empty, since $Y'$ maps to an open subset $\pi(Y')$ of $Q$, whose preimage $\pi^{-1}(\pi(Y'))$ has to intersect $\mathscr{U}$ which is dense in $X$. Better still, by the same reasoning, we see that $\pi_2(\Gamma')$ contains an open $G$-invariant subset $\mathscr{U}'$ of $\mathscr{U}$. Since $\pi_2(\Gamma')$ is constructible by Chevalley's Theorem, it contains a $G$-invariant Zariski-open subset $\mathscr{U}'$ of $\mathscr{U}$. Therefore, we can assume without loss of generality that $(\pi_2)|_{\Gamma'}$ is surjective.

Let $\overline{Y'}$ be a projective normal $G$-equivariant completion of $Y'$ and let $\overline{\Gamma'}$ be the closure of $\Gamma'$ in $\overline{Y'}\times \mathscr{U}$. Similar to $(3)$ above, we can assume that for every $u \in \mathscr{U}$ there exists some $y \in Y'\cap \overline{G\acts u}$ such that
\begin{equation}\label{eq:closureoforbit}
(\pi_2|_{\overline {\Gamma'}})^{-1}(u) =  \overline{(\pi_2|_{\Gamma'})^{-1}(u) } = \overline{G\acts y} \times \{u\}.\tag{$\diamond$}
\end{equation}

Since $\overline{\Gamma'}$ is irreducible, since $\mathscr{U}$ is smooth, and since there exists a constant $k \in \N$ such that $\dim G\acts y$ equals $k$ for all $y \in Y$, Theorem~\ref{thm:algebraicfamily} implies that the graph $\overline{\Gamma'}$ defines a well-defined algebraic family of compact $k$-cycles in $\overline{Y'}$ para\-me\-trised by $\mathscr{U}$. As a consequence we obtain a $G$-invariant regular map $\psi: \mathscr{U} \to \mathscr{C}_{k,d}(\overline{Y'})$ to the Chow variety $\mathscr{C}_{k,d}(\overline{Y'})$ of cycles of degree $d$ and dimension $k$ in $\overline{Y'}$.

Let $U_R$ be a Rosenlicht subset of $\overline{Y'}$ contained in $Y'$ constructed as in Proposition~\ref{Chow}. Clearly $U_R$ is also a Rosenlicht subset of $Y$. Note that due to the assumption $X=\pi^{-1}(S)$, the quotient $X\hq G$ can no longer assumed to be complete. However, by considering the closure $Y$ of $Y'$ in $X$, we may apply Proposition~\ref{prop:rosenlichtandsemistable} and, after shrinking further if necessary, assume without loss of generality that $U_R/G$ has an open embedding $\imath$ into $Q$.

Next we claim that $\psi(\mathscr{U})$ contains an open subset of $U_R/G \subset \mathscr{C}_{k,d}(\overline{Y'})$. Owing to \eqref{eq:closureoforbit}, it suffices to show that $\pi_1(\Gamma)$ contains a $G$-invariant Zariski-open subset of $Y'$. Aiming for a contradiction, suppose that $\overline{\pi_1 (\Gamma)}$ is a proper $G$-invariant subvariety of $Y'$.
Its image $E \definiere \pi(\overline{\pi_1(\Gamma)})$ under $\pi$ is a closed analytic subset of the Zariski-open subset $\pi(Y')$ of $Q$. Set $\tilde{\mathscr{U}} \definiere \pi^{-1}\bigl( (Y'/G) \setminus E\bigr)$. Since $\mathscr{U}$ is dense in $X$, we have $\mathscr{U} \cap \tilde{\mathscr{U}} \neq \emptyset$. However, this contradicts the definition of $E$.

We can thus assume that $\psi$ maps $\mathscr{U}$ into $U_R/G$. Furthermore, again owing to \eqref{eq:closureoforbit}, we have $\pi|_{\mathscr{U}}= \imath \circ \psi|_{\mathscr{U}}$.
Consequently, the restriction of $\pi$ to $\mathscr{U}$ is algebraic and hence extends to a rational map $\Pi: X \dasharrow Q$. Since $\Pi|_U = \pi$ is a priori holomorphic, we conclude that it is regular. This completes the proof of Theorem~\ref{thm:algebraicmap}.
\end{proof}
The following example shows that the completeness assumption is indeed necessary:
\begin{ex}\label{ex:nonalgebraicautomo}
Let $X = U = Q = \C^2$ and let $\varphi$ be a non-algebraic, holomorphic automorphism of $\C^2$, cf.\ Example~\ref{ex:strangeopenset}. Then, $\varphi: U \to Q$ is an analytic Hilbert quotient for the trivial group action.
\end{ex}
\subsection{Categorical quotients}\label{subsect:categorical}
A $G$-invariant morphism $\pi:X \to Y$ from a $G$-variety $X$ to a variety $Y$ is called a \emph{categorical quotient} if for every $G$-invariant morphism $\phi: X \to Z$ to an algebraic variety $Z$, there exists a uniquely determined morphism $\bar \phi : Y \to Z$ such that $\phi = \bar \phi \circ \pi$.
\begin{lemma}\label{lem:categorical}
Let $G$ be a complex reductive Lie group. Let $X$ be a $G$-irreducible algebraic $G$-variety and $U \subset X$ a $G$-invariant Zariski-open subset such that the analytic Hilbert quotient $\pi: U \to U\hq G$ exists. Assume that $U\hq G$ is a complete algebraic variety. Then we have $(\pi_* \mathscr{O}_U)^G = \mathscr{O}_{U/\negthickspace / G}$, and $\pi$ is a categorical quotient.
\end{lemma}
\begin{proof}
First we note that by the previously proven results the quotient map $\pi$ is morphism of algebraic varieties.
Let $V \subset U\hq G$ be Zariski-open. Then $\pi^{-1}(V)$ is Zariski-open and dense in $U$. The pull-back of any regular function on $V$ yields a $G$-invariant regular function on $\pi^{-1}(V)$. Conversely, let $f \in \mathscr{O}_X(\pi^{-1}(V))^G$. Then $f$ extends to a $G$-invariant rational function $F \in \C(U)^G$. Now Proposition~\ref{prop:meromorphicdescent} implies that the holomorphic function $\bar f \in \mathscr{H}_{U/\negthickspace/ G}(V)$ with $\pi^*(\bar f)=f$ extends to a meromorphic, hence rational function on $U\hq G$ (as a consequence of GAGA, \cite[Thm.\ XII.4.4]{SGA1}). Consequently, $\bar f$ is regular.

For the second statement of the lemma we need to show that for every $G$-invariant regular map $\phi: U \to Z$ to an algebraic variety $Z$ there exists a regular map $\bar \phi: U\hq G \to Z$ with
\begin{equation}\label{eq:categorical}
\phi = \bar\phi \circ \pi.\tag{$\star$}
\end{equation}

So let $\phi: U \to Z$ be given. Without loss of generality, we can assume that $Z$ is complete. Since $\pi: U \to U\hq G$ is an analytic Hilbert quotient, there exists a holomorphic map $\bar \phi: U\hq G \to Z$ such that \eqref{eq:categorical} is fulfilled. However, every holomorphic map between complete algebraic varieties is algebraic, see \cite[Cor.\ XII.4.5]{SGA1}.
\end{proof}
In the next section we prove a more general version of this result for the case of non-reduced schemes, cf.\ the proof of Proposition~\ref{prop:coherenceforsubschemes}.
\section{Coherence of sheaves of invariants}\label{sect:coherence}
The aim of this section is to show the following result:
\begin{thm}[Coherence Theorem]\label{thm:coherence}
Let $G$ be a connected complex reductive Lie group. Let $X$ be an irreducible algebraic $G$-variety with analytic Hilbert quotient $\pi^h: X^h \to Q^h \definiere X^h\hq G$ where $Q$ is a complete algebraic variety. If $\mathscr{F}$ is any coherent algebraic $G$-sheaf on $X$, then $(\pi_* \mathscr{F})^G$ is a coherent algebraic sheaf on $Q$.
\end{thm}
Note that if $\pi: X\to Q$ is actually an algebraic Hilbert quotient, the statement follows from the fact that $\pi$ is by definition an affine map, cf. Section 2 of \cite{HeinznerCoherent}. Furthermore, recall that under the assumptions listed above $\pi$ is a morphism of algebraic varieties (Theorem~\ref{thm:algebraicmap}).
\subsection{Coherence for structure sheaves of closed subschemes}
We will start by considering the following special case of the Coherence Theorem.
\begin{prop}\label{prop:coherenceforsubschemes}
Let $G$, $X$, and $\pi:X\to Q$ be as in Theorem~\ref{thm:coherence}. Furthermore, let $(Z, \mathscr{O}_Z) \hookrightarrow (X, \mathscr{O}_X)$ be a closed $G$-invariant subscheme of $X$. Then, $(\pi_*\mathscr{O}_Z)^G$ is a coherent algebraic sheaf on $Q$.
\end{prop}
\begin{rem}
In fact, as the proof will show, $(\pi_*\mathscr{O}_Z)^G$ is the structure sheaf of the closed subscheme $\pi(Z) \hookrightarrow Q$.
\end{rem}
\begin{proof}[Proof of Proposition~\ref{prop:coherenceforsubschemes}]
The closed subscheme $Z \hookrightarrow X$ is given by a $G$-invariant ideal sheaf $\mathscr{I}_Z$ in $\mathscr{O}_X$. The associated complex space $Z^h$ is the closed $G$-invariant complex subspace of $X^h$ given by $\mathscr{I}_Z^h$. If we endow its image $\pi(Z)$ in $Q^h$ with the canonical structure sheaf $\mathscr{H}_{\pi(Z)}$ induced by the ideal sheaf $(\pi_*^h \mathscr{I}_Z^h)^G$, then the restriction $\pi|_{Z^h}: (Z^h, \mathscr{H}_{Z^h}) \to (\pi(Z), \mathscr{H}_{\pi(Z)})$ is an analytic Hilbert quotient, cf.\ Section~\ref{subsect:Actionsandquotients}. Since $Q$ is a complete algebraic variety by assumption, $(\pi(Z), \mathscr{H}_{\pi(Z)})$ is the complex space associated with the closed subscheme of $Q$ that is given by the uniquely determined ideal subsheaf $\mathscr{I}_{\pi(Z)}$ of $\mathscr{O}_Q$ with $\mathscr{I}_{\pi(Z)}^h = (\pi_*^h \mathscr{I}_Z^h)^G$, see \cite[Thm.\ XII.4.4]{SGA1}. Since $\pi^*(\mathscr{I}_{\pi(Z)}) \subset \mathscr{I}_Z$, the quotient map $\pi$ restricts to a morphism of algebraic schemes $\pi|_Z: (Z, \mathscr{O}_Z) \to (\pi(Z), \mathscr{O}_{\pi(Z)})$.

We claim that $((\pi|_Z)_*\mathscr{O}_Z)^G = \mathscr{O}_{\pi(Z)}$, from which the coherence of $(\pi_*\mathscr{O}_Z)^G$ follows immediatly. Abusing notation, we denote $\pi|_Z$ as $\pi$ in the following argument.

To show the claim, let $U$ be a Zariski-open subset of $\pi(Z)$. Clearly, since $\pi$ is a morphism of algebraic schemes, the pullback $\pi^*(f)$ of a regular function $f \in H^0(U, \mathscr{O}_{\pi(Z)})$ is a $G$-invariant regular function on $\pi^{-1}(U)$.

Conversely, let $f \in H^0(\pi^{-1}(U), \mathscr{O}_Z)^G$. The scheme-theoretic closure $Y = \overline{\pi^{-1}(U)}$ is a closed $G$-invariant subscheme of $Z$, its image $\pi(Y)$, endowed with the canonical structure, is a closed subscheme of $Z \hookrightarrow Q$ (again due to GAGA, \cite[Thm.\ XII.4.4]{SGA1}). Note that $(Y \setminus \pi^{-1}(U))^h$ is analytically rare and $\pi$-saturated in $Y^h$ and that $\pi$ restricts to a regular morphism from $Y$ to $\pi(Y)$. Furthermore, $(U, \mathscr{O}_{\pi(Z)}|_U)$ and $(U, \mathscr{O}_{\pi(Y)}|_U)$ are by construction canonically isomorphic as algebraic schemes. The $G$-invariant regular function $f$ extends to a $G$-invariant rational, hence meromorphic function $F \in H^0(Y^h, \mathscr{M}_{Y^h})^G$. By Proposition~\ref{prop:meromorphicdescent}, there exists a meromorphic function $\hat F$ on $\pi(Y)^h$, such that $(\pi^h)^*(\hat F) = F$. Since $\pi(Y)^h$ is a closed complex subspace of $Q^h$, it is also complete algebraic. As a consequence of GAGA, the meromorphic dunction $\hat F$ is (induced from) a rational function on $\pi(Y)$. Furthermore, its restriction $\hat f \definiere \hat F |_U$ to $U$ is holomorphic, hence regular, i.e., $\hat f \in H^0(U, \mathscr{O}_{\pi(Y)}|_U) = H^0(U, \mathscr{O}_{\pi(Z)}|_U)$, and fulfills $\pi^* (\hat f) = f$.
\end{proof}
\subsection{Coherence for general $G$-sheaves}
In the following proof of the Coherence Theorem we use a standard procedure (see e.g.\ \cite[Lemma 3.1]{Japaner}) to reduce the general case of arbitrary coherent algebraic $G$-sheaves to the special case of structure sheaves of $G$-invariant subschemes already considered in the previous section.
\begin{proof}[Proof of Theorem~\ref{thm:coherence}]
Aiming for a contradiction, suppose that there exists a coherent algebraic $G$-sheaf $\mathscr{F}$ on $X$ such that $(\pi_*\mathscr{F})^G$ is not coherent. Then, there exists an affine open subset $Q'$ of $Q$ such that $(\pi_*\mathscr{F})^G|_{Q'} = ((\pi|_{\pi^{-1}(Q')})_*(\mathscr{F}|_{\pi^{-1}(Q')}))^G$ is not coherent. Set $X' \definiere \pi^{-1}(Q')$, $\pi' \definiere \pi|_{X'}$, and $\mathscr{F}' \definiere \mathscr{F}|_{X'}$. There exists a coherent algebraic $G$-subsheaf $\mathscr{G}'$ of $\mathscr{F}'$ with the following properties:
\begin{enumerate}
\item $(\pi'_*\mathscr{G}')^G$ is coherent on $Q'$,
\item the only coherent algebraic $G$-subsheaf $\mathscr{E}$ of $\mathscr{F}'/\mathscr{G}'$ such that $(\pi_*'\mathscr{E})^G$ is coherent is the trivial subsheaf $\mathscr{E}=0$.
\end{enumerate}
 As $\pi_*' (\cdot) ^G$ is left-exact, we may assume that $\mathscr{G}' = 0$ in the following, replacing $\mathscr{F}'$ by $\mathscr{F}'/\mathscr{G}'$, if necessary.

As $(\pi_*'\mathscr{F}')^G$ is not coherent, it is not the zero sheaf. Since $Q'$ is affine, it follows that there exists an element $a \in H^0(X', \mathscr{F}')^G \setminus \{0\}$. The section $a$ yields a $G$-equivariant sheaf morphism $\mathfrak{a}: \mathscr{O}_{X'} \to \mathscr{F}'$. Its image $\mathscr{E} \definiere \mathfrak{a}(\mathscr{O}_{X'})$ is a non-trivial coherent algebraic $G$-subsheaf of $\mathscr{F}'$. By (2) above we conclude that
\begin{equation}\label{eq:notcoherent}
\text{$(\pi'_*\mathscr{E})^G$ is not coherent. }\tag{$\dagger$}
\end{equation}
If $\mathscr{J}$ denotes the kernel of $\mathfrak{a}$, and $Z$ is the $G$-stable closed subscheme of $X'$ defined by $\mathscr{J}$, then $\mathscr{E} \cong \mathscr{O}_{X'}/\mathscr{J} = \mathscr{O}_Z$. With this notation, we obtain the following commutative diagram:
\[\begin{xymatrix}{
Z \ar@{^{(}->}^{\text{closed}}[rr]  & & X'\ar@{^{(}->}^{\text{open}}[rr]\ar^{\pi'}[d] & & X \ar^{\pi}[d] \\
 & & Q' \ar@{^{(}->}^{\text{open}}[rr]& & Q.
}
\end{xymatrix}
\]
Let $\overline{Z}$ be the scheme-theoretic closure of $Z$ in $X$. Since $|\overline Z  \setminus Z|$ does not intersect $X'$, we have $(\pi_*\mathscr{O}_{\overline{Z}})^G|_{Q'} = (\pi'_* \mathscr{O}_Z)^G$. Consequently, Proposition~\ref{prop:coherenceforsubschemes} implies that $(\pi_*'\mathscr{O}_Z)^G \cong (\pi_*'\mathscr{E})^G$ is a coherent algebraic sheaf on $Q'$. This contradicts \eqref{eq:notcoherent}.
\end{proof}
\section{Proof of the main theorem}\label{sect:maintheoremproof}
In this section, we prove Theorem~\ref{thm:mainthm}, after a discussion of the local structure of the quotient map near generic closed orbits.

We recall the setup: let $G$ be a complex reductive Lie group, let $X$ be a $G$-irreducible normal algebraic $G$-variety, and let $U\subset X$ be a $G$-invariant analytically Zariski-open subset of $X$ such that the analytic Hilbert quotient $\pi: U \to U\hq G$ exists, and such that $Q \definiere U\hq G$ is a projective algebraic variety. Then, the Openess Theorem, Theorem~\ref{thm:Zariskiopen}, implies that $U$ is Zariski-open in $X$, and we conclude from the Algebraicity Theorem, Theorem~\ref{thm:algebraicmap}, and from Lemma~\ref{lem:categorical} that the quotient map $\pi$ is a categorical quotient. It remains to show that $\pi$ is in fact an affine map.
\begin{rem}\label{rem:mumfordexample}
It follows from a classical result of Snow \cite[Cor.\ 5.6]{Snow} (and also from \cite[Thm.\ 3]{PaHq}) that in the situation at hand any fibre of $\pi$ naturally carries the structure of an affine algebraic $G$-variety. Examples show that this is not sufficient to guarantee that $\pi$ is an affine map. In fact, Mumford \cite[Ex.\ 0.4]{MumfordGIT} has constructed an example of an algebraic $Sl_2(\C)$-variety $X$ with a geometric categorical quotient $\pi: X \to Q$ such that all stabiliser groups of points in $X$ are trivial, but $\pi$ is not a good quotient; see also \cite[Ex.\ 7.1.1]{BBSurvey}.
\end{rem}
\subsection{An \'etale slice theorem at generic orbits}\label{subsect:etaleslice}
In this section we take a closer look at the local structure of the quotient $\pi: U \to Q$ near generic closed orbits. The main result is:
\begin{prop}[\'Etale slice theorem at generic orbits]\label{prop:Lunaslice}
Let $G$, $X$, $\pi: U \to Q$ be as in Theorem~\ref{thm:mainthm}. Let $p \in \pi^{-1}(Q_{\mathrm{princ}})$ be a point whose orbit $G\acts p$ is closed in the maximal Luna stratum $\pi^{-1}(Q_{\mathrm{princ}})$ of $U$. Then, there exists an algebraic $G$-variety T with algebraic Hilbert quotient $T \to T\hq G$ and a regular $G$-equivariant map $\varphi:T \to U$ onto a $\pi$-saturated Zariski-open neighbourhood $W$ of $p$ in $\pi^{-1}(Q_{\mathrm{princ}})$ such that every point $q \in T$ has a saturated open neighbourhood $V$ in $T$ with the property that $\varphi|_V: V \to \varphi(V)$ is a biholomorphism onto the $\pi$-saturated open subset $\varphi (V)$ of $\pi^{-1}(Q_{\mathrm{princ}})$. In particular, in the commutative diagram
\[\begin{xymatrix}{
T \ar^{\varphi}[r] \ar_{\pi_T}[d]& W \ar^{\pi}[d] \\
T\hq G \ar^{\overline{\varphi}}[r]& W \hq G,
}
\end{xymatrix}
\]
both $\varphi$ and $\overline{\varphi}$ are \'etale. Moreover, $T$ is $G$-equivariantly isomorphic to $W \times_{W /\negthickspace/ G} (T \hq G)$.
\end{prop}
\begin{proof}
Since $H$ is reductive, there exists a Zariski-locally closed $H$-invariant subvariety $S$ of $U$ such that the natural map $\varphi: G\times_H S \to U$ is locally biholomorphic at the point $[e,p] \in G \times_H S$, cf.\ the proof of Theorem 4.6 in \cite{PaHq}. Furthermore, shrinking $S$, we may assume that $S$ is affine, and that $\varphi (G\times_H S)$ is $\pi$-saturated in $U$. As a consequence, $ T \definiere G\times_H S$ is an affine $G$-variety and hence, the algebraic Hilbert quotient $T \to  T\hq G \cong S\hq H$ exists. In the following a \emph{saturated} subset of $T$ always refers to a subset that is saturated with respect to this algebraic Hilbert quotient.

The set $E \definiere \{q \in T \mid \varphi \text{ is not \'etale at }q \}$ is a closed $G$-invariant subvariety of $T$. First substituting $T$ by $ T \setminus \mathcal{S}_G(E)$, and then shrinking further, we may assume that $\varphi$ is \'etale at every point of $T$, and, shrinking $U$, that $\varphi$ is surjective onto $U$.

Applying the holomorphic Slice Theorem \cite[Section 6]{HeinznerGIT} at $p$ yields that the principal orbit type of closed $G$-orbits in $T$ coincides with the principal orbit type of closed orbits in $U$. Since the corresponding Luna stratum $\pi^{-1}(Q_{\mathrm{princ}})$ in $U$ is Zariski-open, again we can assume that all the isotropy groups of closed orbits in $T$ and $U$ are conjugate to $H$ in $G$. This yields the subset $W$ of $U$ with the desired properties.

The image $\varphi(G\acts q)$ of any closed orbit $G\acts q$ in $T= \varphi^{-1}(W)$ is closed in $U$. Indeed, aiming for a contradiction, suppose that there exists a $y \in U$ such that $G \acts y \subsetneq \overline{G \acts \varphi (q)}$ is closed. Then, $\dim (G\acts y) < \dim (G \acts \varphi(q)) \leq \dim (G\acts q)$. However, $\dim G\acts y = G \acts q$, since $G_y$ and $G_q$ are both conjugate to $H$ in $G$. This is a contradiction.

Again due to the fact that all stabiliser groups of closed orbits are conjugate, for every closed orbit $q \in T$, the restriction $\varphi|_{G\acts q}: G\acts q \to G \acts \varphi(q)$ is an isomorphism onto the closed orbit $G \acts \varphi(q) \subset U$. We may thus apply the slice theorem at any point $q \in T$ with closed orbit $G \acts q$ to obtain the open neighbourhood $V$ of $q$ with the desired properties.

Set $Z \definiere W \times_Q (T \hq G)$ and consider the regular map $\chi: T \to Z$, $t \mapsto (\varphi(t), \pi_T(t))$. Since locally over the quotient $\varphi$ is an isomomorphism, $\chi$ is bijective. Let $V \subset T\hq G$ be an open subset such that $\pi_T^{-1}(V)$ is $G$-equivariantly isomorphic to the $\pi$-saturated open subset $\widetilde{V} \definiere \varphi(\pi_T^{-1}(V)) \subset U$ via $\varphi$. We denote the inverse of $\varphi|_{\pi_T^{-1}(V)}$ by $\psi$. Let $\widetilde Z \definiere \widetilde V \times_{\overline{\varphi}(V)} V$. It follows that $\chi|_{\pi_T^{-1}(V)}: {\pi_T^{-1}(V)} \to \widetilde Z$ is biholomorphic with holomorphic inverse given by $(x, [t]) \mapsto \psi(x)$. We have thus shown that $\chi$ is a bijective, everywhere locally biholomorphic, regular map and hence biregular.
\end{proof}
Luna's slice theorem \cite{LunaSlice} applied to the quotient $\pi_T: T \to T \hq G$ together with Proposition~\ref{prop:Lunaslice} now yields the following result:
\begin{cor}\label{cor:quotientfibrebundle}
Let $G$, $X$, $\pi:U\to Q$ be as in Theorem~\ref{thm:mainthm}. Then, the restriction of the quotient map to the maximal Luna stratum
\[\pi|_{\pi^{-1}(Q_{\mathrm{princ}})}: \pi^{-1}(Q_{\mathrm{princ}}) \to Q_{\mathrm{princ}}\]
realises $\pi^{-1}(Q_{\mathrm{princ}})$ as an algebraic \'etale locally trivial fibre bundle over $Q_{\mathrm{princ}}$ with fibre the affine algebraic $G$-variety $G \times_H (\pi_S^{-1}(\pi_S^{\ }(p))$. Here, $p \in \pi^{-1}(Q_{\mathrm{princ}})$ is any point with closed orbit $G\acts p$ in $\pi^{-1}(Q_{\mathrm{princ}})$, $S$ is an affine local slice at $p$, and $\pi_S: S \to S \hq H$ denotes the good quotient for the $H$-action on $S$.
\end{cor}
\subsection{Proof of Theorem~\ref{thm:mainthm}}\label{subsect:maintheoremproof}
First we recall the following well-known corollary of Zariski's Main Theorem:
\begin{lemma}\label{lem:etalefinite}
Let $e: X \to Y$ be a quasi-finite morphism between irreducible algebraic varieties. Then, there exists a non-empty Zariski-open subset $W \subset Y$ such that $e|_{e^{-1}(W)}: e^{-1}(W) \to W$ is a finite map.
\end{lemma}
We now apply this to our setup:
\begin{lemma}\label{lem:locallygood}
Let $G$, $X$, $\pi:U\to Q$ be as in Theorem~\ref{thm:mainthm}. Then, there exists a non-empty Zariski-open subset $W \subset Q$ such that $\pi|_{\pi^{-1}(W)}: \pi^{-1}(W) \to W$ is an affine map.
\end{lemma}
\begin{proof}
Since we are only interested in the behaviour of $\pi$ over a big open set of $Q$, by Corollary~\ref{cor:quotientfibrebundle} we may assume that there exist an irreducible variety $B$, an affine algebraic $G$-variety $F$ with $\C[F]^G = \C$, and an \'etale $G$-equivariant surjective map $\psi: B \times F \to U$, such that the induced map $\bar{\psi}: B \to Q$ is \'etale, and such that in the commutative diagram
\[\begin{xymatrix}{
B \times F \ar^{\psi}[r]\ar[d]  & U \ar^{\pi}[d] \\
B \ar^{\bar{\psi}}[r]    &   Q,
}
\end{xymatrix}
\]
the trivial product $B\times F$ is $G$-equivariantly isomorphic to the fibre product $U \times_Q B$ over $B$. By Lemma~\ref{lem:etalefinite}, there exists a non-empty Zariski-open subset $W\subset Q$ such that $\bar{\psi}|_{\bar \psi ^{-1}(W)}: \bar \psi ^{-1}(W) \to W$ is a finite \'etale covering, and hence proper. Arguing exactly as in the proof of Theorem 8.6 in \cite{PaHq}, we conclude that $\psi|_{\bar \psi ^{-1}(W) \times F}: \bar \psi ^{-1}(W) \times F \to \pi^{-1}(W)$ is likewise proper and therefore finite.

Let now $A \subset W$ be an affine open subset. Then, since $F$ is affine and $\bar{\psi}$ is finite,
$\psi|_{\bar{\psi}^{-1}(A) \times F}: \bar{\psi}^{-1}(A) \times F \to \pi^{-1}(A)$
is a finite map from an affine variety onto $\pi^{-1}(A)$. We hence infer from Chevalley's Theorem \cite[III.Ex.\ 4.2]{Hartshorne} that the preimage $\pi^{-1}(A)$ is also affine.
\end{proof}
The following lemma will be used in the subsequent proof of Theorem~\ref{thm:mainthm}. It is proven exactly as Lemma~\ref{lem:basicconstruction}.
\begin{lemma}\label{lem:sectionstomaps}
Let $X$ be an algebraic $G$-variety with categorical quotient $\pi: X \to Q$. Let $V$ be a finite-dimensional $G$-module, let $\mathscr{L}$ be the sheaf of sections of a line bundle $L$ on $Q$, and let $n\in \N$. Then, every section $s \in H^0\bigl(Q, \pi_*(V \otimes \pi^*\mathscr{L}^{\otimes n})^G\bigr)$ yields a $G$-equivariant algebraic map $\sigma: X \to V \otimes L^{\otimes n}$.
\end{lemma}
\begin{proof}[Proof of Theorem~\ref{thm:mainthm}]
In view of the discussion in Section~\ref{sect:embeddingGspaces}, it suffices to show that for a fi\-nite-\-dimensional $G$-module $V$ every section $s \in H^0\bigl(Q^h, \pi_*^h(V \otimes (\pi^h)^*(\mathscr{L}^h)^{\otimes n})^G\bigr)$ yields a $G$-equivariant \emph{algebraic} map $\sigma: U \to V \otimes L^{\otimes n}$.

The sheaf $\mathscr{F} \definiere V \otimes \pi^*\mathscr{L}^{\otimes n}$ is a coherent algebraic $G$-sheaf on $U$. By the Coherence Theorem, Theorem~\ref{thm:coherence}, $(\pi_*\mathscr{F})^G$ is a coherent algebraic sheaf on $Q$. The analytic sheaf $\mathscr{F}^h = V \otimes (\pi^h)^*(\mathscr{L}^h)^{\otimes n}$ associated with $\mathscr{F}$ is a coherent analytic $G$-sheaf on $U^h$, and consequently,  $(\pi_*^h\mathscr{F}^h)^G$ is a coherent analytic sheaf on $Q^h$ by the Coherence Theorem of \cite{HeinznerCoherent}. Since $Q$ is projective algebraic, GAGA \cite[Thm.\ 3]{GAGA} asserts that there exists a coherent algebraic sheaf $\mathscr{G}$ on $Q$ such that $(\pi_*^h\mathscr{F}^h)^G = \mathscr{G}^h$.

It follows from a result of Neeman \cite{NeemanWeakGAGA} (see also Lemma 2 of \cite{NeemanAnalytic}) that the natural morphism $\nu: \bigl((\pi_*\mathscr{F})^G \bigr)^h \to (\pi_*^h\mathscr{F}^h)^G = \mathscr{G}^h$ is injective. Again by GAGA, the morphism $\nu$ is induced by an injective morphism $(\pi_*\mathscr{F})^G \hookrightarrow \mathscr{G}$ of coherent algebraic sheaves, i.e., $(\pi_*\mathscr{F})^G$ is an algebraic subsheaf of $\mathscr{G}$.

By Lemma~\ref{lem:locallygood}, there exists a Zariski-open subset $W$ of $Q$ such that $\pi|_{\pi^{-1}(W)}: \pi^{-1}(W) \to W$ is an algebraic Hilbert quotient. From the GAGA Theorem for quotient morphisms \cite[Thm.\ 8]{NeemanAnalytic} we deduce that
\begin{equation}\label{eq:okonopen}
\bigl((\pi_*\mathscr{F})^G\bigr)^h|_{W^h}= (\pi_*^h\mathscr{F}^h)^G|_{W^h}=\mathscr{G}^h|_{W^h}. \tag{$\ddagger$}
\end{equation}
Next, we consider the exact sequence $0 \to (\pi_*\mathscr{F})^G \to \mathscr{G} \to \mathscr{G}/(\pi_*\mathscr{F})^G \to 0$.
It follows from exactness of the $(\cdot)^h$-functor and from the equality \eqref{eq:okonopen} that $\bigl(\mathscr{G}/(\pi_*\mathscr{F})^G\bigr)^h|_{W^h} = \mathscr{G}^h/\bigl((\pi_*\mathscr{F})^G\bigr)^h|_{W^h} = 0$.  By faithfulness of the $(\cdot)^h$-functor we conclude that
\begin{equation}\label{eq:sheavesonopen}
\mathscr{G}|_W = (\pi_*\mathscr{F})^G|_W \tag{$\sharp$}.
\end{equation}
Let now $s \in H^0 \bigl(Q^h, (\pi_*^h\mathscr{F}^h)^G \bigr) = H^0 \bigl( Q^h, \mathscr{G}^h\bigr)$ and $\sigma: X \to V \otimes L^{\otimes n}$ the corresponding $G$-equivariant holomorphic map, cf.\ Section~\ref{sect:embeddingGspaces}. GAGA says that $s \in H^0 \bigl(Q, \mathscr{G} \bigr)$. Consequently, \eqref{eq:sheavesonopen} implies that $s|_W \in H^0\bigl(W, (\pi_*\mathscr{F})^G \bigr)$. We conclude that the restriction of $\sigma$ to $\pi^{-1}(W)$ is an algebraic map, cf.\ Lemma~\ref{lem:sectionstomaps}. Since $\sigma $ is a priori holomorphic, this implies that it is regular. This concludes the proof of Theorem~\ref{thm:mainthm}.
\end{proof}
\part{Applications}\label{part:applications}
\section{Analytic Hilbert quotients and Geometric Invariant Theory}
In general, it is a fundamental problem of Geometric Invariant Theory to find all open $G$-invariant subsets $U$ of a given algebraic $G$-variety such that a quotient for the induced action on $U$ exists and has certain prescribed properties. The main result of this section, Theorem~\ref{thm:linearisedWeil}, asserts that we do not obtain ''more'' quotients by considering projective analytic Hilbert quotients than by considering linearised Weil divisors in the sense of Hausen's generalisation of Mumford's GIT, which we recall below.
\subsection{GIT based on Weil divisors}\label{subsect:GITwithWeil}
Here, we summarise the definitions and results of \cite{HausenGITwithWeildivisors} that we will use in the following. Let $X$ be a normal $G$-variety and $D$ be a $G$-invariant Weil divisor on $X$. The graded $\mathscr{O}_X$-algebra $\mathcal{A} \definiere \bigoplus_{n \geq 0} \mathscr{O}_X(nD)$ carries a canonical $G$-\emph{linearisation}, i.e., if $\alpha: G \times X \to X$ denotes the action map, and $\mathrm{pr}_X: G \times X \to X$ the canonical projection, then there exists an isomorphism $\phi: \alpha^*\mathcal{A} \to \mathrm{pr}_X^*\mathcal{A}$ of graded $\mathscr{O}_{G \times X}$-algebras such that $\phi$ is the identity in degree zero and such that a certain cocycle condition is fulfilled, see \cite[Def.\ 1.1]{HausenGITwithWeildivisors}.
The group $G$ acts on $\mathcal{A}$ via its usual action on the function field $\C(X)$, i.e., $g\cdot f (x) = f(g^{-1}\acts x)$. For a homogeneous section $f \in \mathcal{A}$ we define the \emph{zero set} to be $Z(f) \definiere \mathrm{supp}\bigl(\mathrm{div}(f) + D\bigr)$. Given a $G$-linearised Weil divisor $D$ as above, a point $x \in X$ is called \emph{semistable} with respect to $D$ if there exists an integer $n>0$ and a $G$-invariant $f \in \mathscr{O}_X(nD)$ such that $X\setminus Z(f)$ is an affine neighbourhood of $x$ on which $D$ is Cartier. The set of all semistable points is denoted by $X(D, G)$, or, in case $D$ is Cartier with corresponding line bundle $L$, by $X(L,G)$. By the following result this concept of semistability yields all open subsets admitting a quasiprojective quotient space:
\begin{thm}[Thm.\ 3.3 of \cite{HausenGITwithWeildivisors}]
Let a reductive group $G$ act on a normal variety $X$.
\begin{enumerate}
\item For a $G$-linearised divisor $D$ on $X$ there exists a good quotient $X(D,G) \to X(D,G) \hq G$ with a quasi-projective variety $X(D,G)\hq G$.
\item If $U \subset X$ is Zariski-open, $G$-invariant, and has a good quotient $U \to U\hq G$ with $U\hq G$ quasiprojective, then $U$ is a $G$-saturated subset of the set $X(D,G)$ of semistable points of some $G$-invariant, canonically $G$-linearised Weil divisor $D$ on $X$.
\end{enumerate}
\end{thm}
\subsection{Analytic Hilbert quotients via $G$-linearised Weil divisors} As a direct corollary of Theorem~\ref{thm:mainthm} and of the results discussed in the previous section, we obtain
\begin{thm}\label{thm:linearisedWeil}
Let $G$ be a complex reductive Lie group, let $X$ be a $G$-irreducible normal algebraic $G$-variety, and let $U \subset X$ be a nonempty $G$-invariant analytically Zariski-open subset of $X$ such that the analytic Hilbert quotient $\pi: U \to U\hq G$ exists. If $U\hq G$ is a projective algebraic variety, then there exists a $G$-linearised Weil divisor $D$ on $X$ such that $U = X(D, G)$.
\end{thm}
\newpage
\section{K\"ahler quotients and Geometric Invariant Theory}\label{sect:KaehlerquotientsGIT}
\subsection{Momentum maps and analytic Hilbert quotients}\label{subsect:intromomentummapquotients}
Let $K$ be a Lie group with Lie algebra $\mathfrak{k}$. Let $X$ be a reduced complex $K$-space endowed with a smooth $K$-invariant K\"ahler structure $\omega = \{\rho_j\}_{j}$ in the sense of Grauert, cf. \cite{Extensionofsymplectic}, \cite{PaHq}. A \emph{momentum map} with respect to $\omega$ is a $K$-equivariant smooth map $\mu: X \to \mathfrak{k}^*$ such that
\[d\mu^\xi = \iota_{\xi_Y}\omega_Y\]
holds for every $K$-stable complex
submanifold $Y$ of $X$ and for every $\xi \in \mathfrak{k}$. Here,  $\iota_{\xi_Y}$ denotes contraction
with the vector field $\xi_Y$ on $Y$ that is induced by the $K$-action, $\omega_Y$ denotes the K\"ahler form induced on $Y$, and the function $\mu^\xi: X
\to \R$ is given by $\mu^\xi(x)=\mu(x)(\xi)$. We call the action of $K$ on a complex $K$-space with
$K$-invariant K\"ahler structure $\omega$ \emph{Hamiltonian} if it admits a
momentum map with respect to $\omega$.

Let $K$ now denote a compact Lie group and $G=K^\C$ its complexification. Then $G$ is reductive. Let $X$ be a reduced holomorphic $G$-space with $K$-invariant K\"ahler structure $\omega$. We call $X$ a
\emph{Hamiltonian $G$-space}, if the $K$-action is Hamilto\-ni\-an with respect to $\omega$. Given a
Hamiltonian $G$-space with momentum map $\mu: X \to \mathfrak{k}^*$, set
\[X(\mu) \definiere \{ x \in X \mid \overline{G\acts x} \cap \mu^{-1}(0) \neq \emptyset\}.\]We
call $X(\mu)$ the set of \emph{semistable points} with respect to $\mu$ and the $G$-action. We collect the main results about Hamiltonian $G$-spaces in
\begin{thm}[\cite{ReductionOfHamiltonianSpaces},
\cite{Extensionofsymplectic}]\label{propertiesmomentumquotients}
Let $X$ be a Hamiltonian $G$-space. Then
\begin{enumerate}
\item the set $X(\mu)$ of semistable points is open and $G$-invariant, and the analytic Hilbert quotient $\pi: X(\mu) \to X(\mu)\hq G$ exists,
\item the inclusion $\mu^{-1}(0) \hookrightarrow X(\mu)$ induces a homeomorphism $\mu^{-1}(0)/K \simeq X(\mu)\hq G$,
\item the complex space $X(\mu)\hq G$ carries a K\"ahler structure that is smooth along a natural stratification of $X(\mu)\hq G$.
\end{enumerate}
\end{thm}
An \emph{algebraic Hamiltonian $G$-variety} is an algebraic $G$-variety $X$ such that the associated holomorphic $G$-space $X^h$ is a Hamiltonian $G$-space. In particular, $X^h$ is assumed to carry a $K$-invariant K\"ahler structure.
\subsection{K\"ahler quotients of algebraic Hamiltonian $G$-varieties}
Using the results obtained above, we can now give a proof of Theorem~\ref{thm:mainthmHamiltonian}. Let us recall the statement:
\begin{thm1.1}[Algebraicity of momentum map quotients]
Let $G=K^\C$ be a complex reductive Lie group and let $X$ be a $G$-irreducible algebraic Hamiltonian $G$-variety with at worst $1$-rational singularities. Assume that the zero fibre $\mu^{-1}(0)$ of the momentum map $\mu: X \to \mathfrak{k}^*$ is nonempty and compact. Then
\begin{enumerate}
\item the analytic Hilbert quotient $X(\mu)\hq G$ is a projective algebraic variety,
\item the set $X(\mu)$ of $\mu$-semistable points is algebraically Zariski-open in $X$,
\item the map $\pi: X(\mu) \to X(\mu)\hq G$ is a good quotient,
\item there exists a $G$-linearised Weil divisor $D$ on $X$ such that $X(\mu)= X(D, G)$.
\end{enumerate}
\end{thm1.1}
\begin{proof}
Part (1) is \cite[Thm.~1]{PaHq}, part (2) is \cite[Thm.~2]{PaHq}. Part (3) follows directly from Theorem \ref{thm:mainthm}, and part (4) is deduced immediatly from Theorem~\ref{thm:linearisedWeil}.
\end{proof}
\begin{rems}\label{rems:finalremarks}
1.) Note that we cannot apply Theorem~\ref{thm:Zariskiopen} to $X(\mu)$, as it is a priori not evident that $X(\mu)$ is an analytically Zariski-open subset of $X$. However, the proof of part (1) in \cite{PaHq} is quite similar to the proof of Theorem~\ref{thm:Zariskiopen} in Section~\ref{subsect:proofofZariskiopen}.

2.) In \cite{MomentumProjectivity} the results summarised in Theorem~\ref{thm:mainthmHamiltonian} were obtained for the special case of smooth projective algebraic Hamiltonian $G$-varieties using Hodge theory and exhaustion properties of strictly plurisubharmonic functions. Furthermore, it was shown that one can choose $D$ to be an ample Cartier divisor in this setup. We will prove a generalisation of this result in the case of semisimple group actions in Section~\ref{subsect:semisimple}.

3.) In \cite{PaHq}, an example of a non-algebraic, non-compact momentum map quotient of a smooth algebraic Hamiltonian $\C^*$-variety was constructed. Hence, there is no extension of the result above to the case of non-compact momentum map quotients.

4.) See \cite{PaHq} or \cite{GrebSingularities} for basic properties and examples of $1$-rational singularities. The class of complex spaces with $1$-rational singularities is the natural class of singular varieties to which projectivity results for K\"ahler Moishezon manifolds generalise, cf.\ \cite{ProjectivityofMoishezon}. Furthermore, it is also natural from the point of view of equivariant geometry, since it is stable under taking good quotients and analytic Hilbert quotients, see \cite{GrebSingularities} and \cite{GrebAnalyticSingularities}.

5.) The assumption on the singularites of $X$ is used in \cite{PaHq} to show that the quotient $X(\mu)\hq G$ is projective applying the results of Namikawa referred to above. It is necessary in order to obtain projectivity: by the example given below there exist (necessarily singular) non-projective, complete K\"ahlerian algebraic varieties. These are momentum map quotients for the trivial group action.
\end{rems}
In the following, we construct such an example of a K\"ahlerian non-projective proper algebraic surface. While there exist quite a few constructions of non-projective proper algebraic varieties and Moishezon spaces, the author was not able to find a reference where the existence of K\"ahler structures on the resulting examples is discussed. In this sense the surface described below, which was first constructed by Schr\"oer, is the first example of a K\"ahlerian non-projective proper algebraic variety.
\begin{ex}\label{ex:nonprojective}
Using the construction described in \cite[Sect.\ 2]{SchroerNonprojective} one finds a non-projective normal proper algebraic surface $X$ having exactly two singular points $x_1, x_2 \in X$ and a resolution $\pi: \widetilde X \to X$ with the following properties:
\begin{itemize}
\item $\widetilde X$ is a projective surface; in particular, $X$ is in Fujiki class $\mathscr{C}$,
\item the exceptional divisors $E_1$ and $E_2$ are smooth curves of genus one,
\item $E_1^2 = E_2^2 = -1$.
\end{itemize}
We claim that $X$ is K\"ahler. According to the main result of \cite{FujikiKaehlerSurfaces} this follows as soon as we can show that the natural map $\varphi: R^1\pi_*\R_{\widetilde X} \to R^1\pi_*\mathscr{O}_{\widetilde X}$ is surjective. Since $\pi$ is proper, since both sheaves are supported on $\{x_1, x_2\}$, since $\dim_\R\bigl((R^1\pi_*\R_{\widetilde X})_{x_j}\bigr) = 2$ and since $\varphi$ is always injective, see \cite[Prop.\ 6]{FujikiKaehlerSurfaces}, it suffices to show that for $j=1,2$ there exists a strictly pseudoconvex neighbourhood $T_j$ of $E_j$ in $\widetilde X$ such that $H^1\bigl(T_j, \mathscr{H}_{T_j} \bigr) = \C$. Note that this will show in particular that the singularities of $X$ are not $1$-rational. Since the problem is local, we will omit the subscript $j$ in the following discussion. We may assume that $X$ is a closed subset of $\C^N$ and $x = 0 \in \C^N$. Intersecting a small ball centered at $0$ with $X$ and pulling back the resulting open (Stein) subset to $\widetilde X$ via $\pi$, we obtain a strictly pseudoconvex neighbourhood $T$ of $E$ in $\widetilde X$. By a fundamental result of Grauert~\cite[§4, Satz 1]{GrauertModifikationen} there exists an $m_0 \in \N^{>0}$ such that for the $m_0$-th infinitesimal neighbourhood $m_0E$ of $E$ in $\widetilde X$ we have $H^1\bigl(T, \mathscr{H}_{T} \bigr) \cong H^1\bigl(m_0E, \mathscr{H}_{m_0E}\bigr)$.  For any $m \in \mathbb{N}^{>0}$ consider the exact sequence
\begin{equation}\label{eq:infneighbourhoods}
0 \to \mathscr{N}^{*\otimes m} \to \mathscr{H}_{(m+1)E} \to \mathscr{H}_{mE} \to 0,\tag{$\spadesuit$}
\end{equation}
where $\mathscr{N}^*$ is dual of (the locally free sheaf associated with) the normal bundle $\mathscr{N}$ of $E$ in $\widetilde X$. Since the degree of $\mathscr{N}$ is negative, we have $H^1\big(E, \mathscr{N}^{*\otimes m} \bigr) \overset{\text{dual}}{\sim} H^0\bigl(E, \mathscr{N}^{\otimes m} \bigr)=\{0\}$ for every $m \in \N^{>0}$. It follows from the long exact sequence of cohomology associated with \eqref{eq:infneighbourhoods} that $H^1\bigl(T, \mathscr{H}_{T} \bigr) \cong H^1\bigl(m_0E, \mathscr{H}_{m_0E}\bigr) \cong  H^1\bigl(E, \mathscr{H}_{E} \bigr) = \C$. Hence, we have shown that $X$ is a K\"ahlerian non-projective proper algebraic surface.
\end{ex}
\subsection{Finiteness of momentum map quotients}
Since the number of Zariski-open subsets of an algebraic $G$-variety that admit a good quotient with complete quotient space is finite by \cite{BBFiniteness}, Theorem~\ref{thm:mainthmHamiltonian} implies
\begin{cor}[Finiteness of momentum map quotients]
Let $G=K^\C$ be a complex reductive Lie group and let $X$ be a $G$-irreducible algebraic $G$-variety with at worst $1$-rational singularities. Then there exist only finitely many (necessarily Zariski-open) subsets of $X$ that can be realised as the set of $\mu$-semistable points with respect to some $K$-invariant K\"ahler structure and some momentum map $\mu: X \to \mathfrak{k}^*$ with compact zero fibre $\mu^{-1}(0)$.
\end{cor}
\subsection{Momentum map quotients modulo semisimple groups}\label{subsect:semisimple}
When studying the momentum geometry of a Hamiltonian $G$-space, especially via the gradient flow of the norm square of the momentum map, it is a common assumption that the restriction of the momentum map to a maximal torus is proper. Note that this is rather restrictive, as the discussion in \cite[sect.\ 2.6]{Doktorarbeit} shows. In this section, we study the implications of the weaker hypothesis that the restricted momentum map has compact zero fibre.

We will frequently use the following notation: if $X$ is a Hamiltonian $K^\C$-space with momentum map $\mu: X \to \mathfrak{k}^*$, and $M < K$ is any closed subgroup with Lie algebra $\mathfrak{m}$, the composition $\mu_M \definiere \imath_{\mathfrak{m}}^* \circ \mu : X \to \mathfrak{m}^*$ of $\mu$ with the map $\imath_{\mathfrak{m}}^*$ dual to the inclusion $\imath_{\mathfrak{m}}: \mathfrak{m} \to \mathfrak{k}$ is a momentum map for the $M$-action. The corresponding set of semistable points $\mathcal{S}_{M^\C}(\mu^{-1}_M(0))$ will be denoted by $X(\mu_M)$.
The main result is
\begin{thm}\label{thm:semisimplecompact}
Let $G$ be a connected semisimple complex Lie group with maximal compact subgroup $K$, and let $X$ be an irreducible algebraic Hamiltonian $G$-variety with at worst $1$-rational singularities and momentum map $\mu: X \to \mathfrak{k}^*$. Let $T$ be a maximal torus of $K$ such that $\mu^{-1}_T(0)$ is non-empty and compact. Then $X$ is projective, we have $X=G\acts X(\mu_T)$, and there exists a $G$-linearised ample line bundle $L$ on $X$ such that $X(\mu_T) = X(L, T^\C)$, and if $X(\mu)\neq \emptyset$, then $X(\mu)$ coincides with $X(L,G)$.
\end{thm}
\begin{proof}
The smallest complex subgroup of $G$ containing $T$ is isomorphic to the complexification $T^\C$ of $T$, and is a maximal algebraic torus of $G$. We denote by $N \definiere N_G(T^\C)$ its normaliser in $G$, and note that $N$ is isomorphic to the complexification of $N_K(T)$. Furthermore, the connected component $N_K(T)^0$ of the identity in $N_K(T)$ equals $T$, and consequently $N^0 = T^\C$. The considerations in the proof of Theorem 1 in Section 5.1 of \cite{PaHq} now show that $X\bigl(\mu_{N_K(T)}\bigr)$ equals $X(\mu_T)$, and consequently, that $X(\mu_T)$ is $N$-stable. By Theorem~\ref{thm:mainthmHamiltonian}, there exists a good quotient $\pi_{T^\C}: X(\mu_T) \to X(\mu_T)\hq T^\C$ with projective quotient variety $X(\mu_T)\hq T^\C$. Since $G$ is semisimple by assumption, it now follows from \cite[Thm.~5.4]{HausenGITwithWeildivisors} that $X(\mu_T)=X(L,T^\C)$ for some $G$-linearised ample line bundle $L$ on $X$, that $X=G\acts X(\mu_T)$, and that $X$ is projective.

Since $\mu^{-1}(0) \subset \mu^{-1}_T(0)$, the quotient $X(\mu)\hq G$ is likewise projective by Theorem~\ref{thm:mainthmHamiltonian}. It follows from \cite[Lem.\ 2.4]{ReductionInStepsBook} and \cite[Thms 4.1 and 5.2]{HausenGITwithWeildivisors}, respectively, that both $X(\mu)$ and $X(L, G) = \bigcap_{g \in G} g \acts X(L, T^\C)$ are $\pi_{T^\C}$-saturated subsets of $X(\mu_T)=X(L, T^\C)$. Let $q: X(L, G) \to X(L, G)\hq G$ denote the quotient map.
\begin{claim}
 $X(\mu)$ is a $q$-saturated subset of $X(L, G)$.
\end{claim}
Assuming the claim we are done since then the compactness of $X(\mu)\hq G$ observed above implies that $X(L, G)\hq G = q(X(\mu))$, and hence that $X(\mu)=X(L, G)$.

It remains to prove the claim. If $M < K$ is a compact subgroup with complexification $H = M^\C < K^\C=G$, and $g \in G$ is any element, the action of $M'\definiere gMg^{-1}$ on $X$ is Hamiltonian with respect to the $M'$-invariant K\"ahler form $(g^{-1})^*\omega$ with momentum map $\mu_{M'} = \Ad^*(g)\circ \mu_M \circ g^{-1}$ and corresponding set of semistable points $X(\mu_{M'}) = \mathcal{S}_{H'}^{}(\mu^{-1}_{M'}(0))$, where $H' = (M')^\C = gHg^{-1}$. With these notations, we have $X(\mu_{gMg^{-1}}) = X(\mu_{M'}) = g \acts X(\mu_M)$. Since $X(\mu)$ is $G$-invariant, we conclude that for all $g \in G$ we have $X(\mu) = X(\mu_{gKg^{-1}}) \subset X(\mu_{gTg^{-1}}) = g\acts X(\mu_T) = g\acts X(L, T^\C)$. It remains to show that $X(\mu)$ is $q$-saturated in $X(L, G)= \bigcap_{g \in G} g \acts X(L, T^\C)$. Suppose on the contrary that there exists $x \in X(\mu)$ such that $\overline{G\acts x} \cap X(\mu) \subsetneq \overline{G\acts x} \cap X(L, G)$. Let $y \in \bigl(\overline{G\acts x} \cap X(L, G)\bigr) \setminus X(\mu)$ such that $G\acts y$ is closed in $X(L,G)$. Then, the Hilbert Lemma (see e.g.\ \cite[Thm.\ 4.2]{HausenGITwithWeildivisors}) asserts that there exists a $g\in G$ such that $\bigl(\overline{T^\C \acts (g \acts x)} \cap X(L, G) \bigr) \cap G\acts y \neq \emptyset$ contradicting the saturatedness of $X(\mu)$ in $X(\mu_T)$.
\end{proof}
As a consequence of Theorem~\ref{thm:semisimplecompact} we are now in the position to extend the main result of \cite{MomentumProjectivity} to singular varieties in the case of semisimple group actions.
\begin{cor}\label{cor:HMgeneralised}
Let $G$ be a connected semisimple complex Lie group with maximal compact subgroup $K$, and let $X$ be an irreducible projective algebraic Hamiltonian $G$-variety with at worst $1$-rational singularities with momentum map $\mu: X \to \mathfrak{k}^*$. Then there exists an ample $G$-linearised line bundle $L$ on $X$ such that $X(\mu)$ coincides with $X(L, G)$.
\end{cor}
\subsection{Momentum geometry versus GIT}
Given an ample $G$-linearised line bundle $L$ on a quasi-projective variety $X$, choosing a compatible Hermitian metric on $L$ gives rise to a momentum map $\mu$ for a maximal compact subgroup $K$ of $G$ such that $X(L, G) = X(\mu)$. Hence, one might hope that every GIT-quotient of a (quasi-projective) algebraic $G$-variety is also a momentum map quotient. However, Bia{\l}ynicki-Birula and {\'S}wi{\polhk{e}}cicka \cite{BBExoticOrbitSpaces} have constructed a smooth projective algebraic $(\C^* \times \C^*)$-variety $X$ together with an invariant Zariski-open subset $U$ such that a good quotient $\pi: U \to U\hq (\C^* \times \C^*)$ exists with projective algebraic quotient variety $U\hq (\C^* \times \C^*)$ and such that there exists no ample linearised line bundle $L$ on $X$ with $X^{ss}(L) = U$; see also \cite[Ex.\ 6.2]{HausenBeyondAmpleCone}. Hence, the results of \cite{MomentumProjectivity} mentioned in Remark~\ref{rems:finalremarks} above imply that there exists no $(S^1 \times S^1)$-invariant K\"ahler form with momentum map $\mu$ on $X$ such that $X(\mu)$ coincides with $U$. It follows that, given an algebraic $G$-variety, one can construct more open subsets with good quotient using Geometric Invariant Theory than using momentum geometry.

Note however that Theorem~\ref{thm:mainthmHamiltonian} above provides a way of constructing interesting K\"ahler structures (i.e., ones that do not arise as curvature forms of ample line bundles) on GIT-quotients of (quasi-projective) $G$-varieties. Furthermore, it indicates that it is interesting to study the variation of the induced K\"ahler structures, similar to the study of variation of GIT-quotients (e.g.\ see \cite{DolgachevHu}, \cite{Thaddeus}), beyond the work of Fujiki \cite{FujikiKaehler}.
\appendix

\section{The principal Luna stratum}\label{appendix:Luna}
Here we prove Lemma~\ref{lem:LunaZopen} in two steps. Let us recall the setup: let $G$ be a connected complex reductive Lie group, let $X$ be an irreducible normal algebraic $G$-variety, and let $U$ be a $G$-invariant Zariski-open subset of $X$ such that the analytic Hilbert quotient $\pi: U \to Q$ exists. Assume that the quotient $Q$ is projective algebraic. Let $S \definiere Q_{\mathrm{princ}}$ be the principal Luna stratum.
\begin{lemma}\label{lem: LunacontainsZopen}
 $\pi^{-1}(S) \subset X(\mu)$ contains a $G$-invariant algebraically Zariski-open subset of $X$.
\end{lemma}
\begin{proof}
Since $\pi_Y: Y^{ss} \to Q$ is algebraic (cf.\ Proposition~\ref{prop:rosenlichtandsemistable}), $A\definiere \pi_Y^{-1}(S)$ is Zariski-open in $Y$. Furthermore, by Lemma~\ref{lem:principalorbitsclosed}, $A$ consists of points whose orbits are closed in $U$. Let $W$ be a Sumihiro neighbourhood of a point $x \in A$ and let $\psi: W \to  \P(V)$ be a $G$-equivariant embedding of $W$ into the projective space associated with a rational $G$-representation $V$. Let $Z$ be the closure of $\psi(W)$ in $\P(V)$.
Given a Rosenlicht subset $U_R$ of $Z$ as in Proposition~\ref{Chow}, Lemma 6.3 of \cite{PaHq} implies that $\mathcal{S}^Z_G(\psi(A \cap W))\cap U_R$ is constructible in
$U_R$. Therefore,
$\mathcal{S}_G^X(A \cap W)\cap W$ contains a $G$-invariant Zariski-open subset $\widetilde U$ of its closure. By construction, $\pi(A \cap W)$ is open in $X(\mu)\hq
G$, and hence, $\pi^{-1}(\pi(A \cap W)) \cap W$ is an open subset of $X$ that is contained in
$\mathcal{S}_G(A \cap W)\cap W$.
\end{proof}
\begin{prop}\label{prop:Lunaalgebraic}
Let $G$ be a connected complex reductive Lie group and $X$ an algebraic $G$-variety. Let $U$ be a $G$-invariant Zariski-open subset of $X$ such that the analytic Hilbert quotient $\pi: U \to Q$ exists. Let $S \definiere Q^{(H)}$. Assume that
\begin{enumerate}
\item $Q$ is projective algebraic,
\item every closed orbit in $U$ has a Sumihiro neighbourhood,
\item there exists a conjugacy class $(H)$ of isotropy groups of closed orbits in $U$ such that $(H) \leq (G_x)$ for all $x \in U$ with closed orbit $G\acts x \subset U$.
\end{enumerate}
Then $\pi^{-1}(S)$ is Zariski-open in $X$.
\end{prop}
\begin{proof}
Without loss of generality, we can assume that $X=U$. Let $X = \bigcup_{j=1}^m X_j$ be the decomposition of $X$ into irreducible components. Then there exists a $j_0 \in \{1, \dots, m\}$ such that $(H)$ is the principal orbit type for the action of $G$ on $X_j$. As a consequence, $\pi^{-1}(S) \cap X_{j_0}$ is analytically Zariski-open in $X_{j_0}$. Furthermore, by Lemma~\ref{lem: LunacontainsZopen}, it contains a $G$-invariant Zariski-open subset $V_{j_0}$ of $X_{j_0}$. Set $\widetilde {X} \definiere  X \setminus \bigl(V_{j_0} \setminus \bigcup_{j\neq j_0} X_j \bigr)$.
Then either $\widetilde{X} = X \setminus \pi^{-1}(S)$, or there exists a closed orbit in $\widetilde X$ with stabiliser $H$. In the first case, we are done, since $\widetilde {X}$ is algebraic in $X$. In the second case, we notice that $\widetilde X^{ss} = \widetilde X$, that $\widetilde X^{ss} \hq G = \pi (\widetilde X) \subset Q$ is projective algebraic, and that the existence of Sumihiro neighbourhoods is inherited by $\widetilde X$. We thus proceed by Noetherian induction.
\end{proof}
\section*{Acknowledgements}
The author wants to thank Peter Heinzner, Christian Miebach, Gerry Schwarz, and Henrik St\"otzel for fruitful and stimulating discussions. Furthermore, the author wants to thank the referee for his helpful and constructive suggestions.
\def\cprime{$'$} \def\polhk#1{\setbox0=\hbox{#1}{\ooalign{\hidewidth
  \lower1.5ex\hbox{`}\hidewidth\crcr\unhbox0}}}
  \def\polhk#1{\setbox0=\hbox{#1}{\ooalign{\hidewidth
  \lower1.5ex\hbox{`}\hidewidth\crcr\unhbox0}}}
\providecommand{\bysame}{\leavevmode\hbox to3em{\hrulefill}\thinspace}
\providecommand{\MR}{\relax\ifhmode\unskip\space\fi MR }
\providecommand{\MRhref}[2]{%
  \href{http://www.ams.org/mathscinet-getitem?mr=#1}{#2}
}
\providecommand{\href}[2]{#2}

\end{document}